\documentclass[10pt,oneside,a4paper,dvipsnames]{amsart}
\usepackage{amsmath,amssymb, amsthm, amsfonts}
\usepackage{hyperref}
\usepackage{cleveref}

\usepackage{mathtools}
\usepackage{thmtools}
\usepackage{etoolbox}
\usepackage[colorinlistoftodos,textwidth=1in]{todonotes}
\usepackage{xparse}
\usepackage{tikz-cd}
\usepackage{graphicx}
\usepackage{fancyhdr}
\usepackage{array}
\usepackage[all,textures]{xy}
\usepackage{enumitem}
\usepackage{caption}
\usepackage{xcolor}

\usepackage[utf8]{inputenc}    
\usepackage[T1]{fontenc}
\usepackage[mathscr]{euscript}
\usepackage{pdfsync}

\usepackage[backend=biber,style=numeric,giveninits,url=false,doi=false,isbn=false,maxbibnames=99]{biblatex}
\hypersetup{hypertexnames = true, bookmarksdepth = 2, bookmarksopen = true, colorlinks, linkcolor = MidnightBlue, citecolor = MidnightBlue, urlcolor = NavyBlue, pdfstartview={XYZ null null 1}}

\mathchardef\mhyphen="2D

\newcommand*{\WidestLabel}{\scriptsize{$(\iota_{\alpha,\beta}^\bbb)_s \circ -$}}% <-- This should be the widest content
\newcommand*{\TallestLabel}{\scriptsize{$S_{\aaa^\beta}^{\bbb^\alpha}$}}% <-- This should be the tallest content

\newlength{\LabelWidth}
\newcommand*{\LabelCong}[1]{\overset{\makebox[\LabelWidth]{\scriptsize{$#1$}\vphantom{\TallestLabel}}}{\cong}}

\let\oldtocsection=\tocsection
\let\oldtocsubsection=\tocsubsection
\newcommand{\tocsectionformat}{}
\expandafter\def\csname toc@1format\endcsname{\tocsectionformat}
\renewcommand{\tocsectionformat}{\bfseries \oldtocsectionformat}
\newcommand{\settocsectionformat}[1]{\renewcommand{\tocsectionformat}{#1}}
\settocsectionformat{\bfseries}
\renewcommand{\tocsection}[2]{\hspace{0em}\bfseries\oldtocsection{#1}{#2}\tocsectionformat}
\renewcommand{\tocsubsection}[2]{\hspace{1.9em}\oldtocsubsection{#1}{#2}}

\DeclareDocumentCommand\issue{g}{\todo[size=\footnotesize,color = green!40]{Issue\IfNoValueF{#1}{: #1}}}
\DeclareDocumentCommand\tobedone{g}{\todo[size=\footnotesize,color = yellow!50]{To be done\IfNoValueF{#1}{: #1}}}
\DeclareDocumentCommand\notationissue{g}{\todo[size=\footnotesize,color = red!30]{Notation?\IfNoValueF{#1}{: #1}}}
\DeclareDocumentCommand\doubt{g}{\todo[size=\footnotesize,color = blue!10]{Doubt\IfNoValueF{#1}{: #1}}}
\DeclareDocumentCommand\observation{g}{\todo[size=\footnotesize,color = orange!10]{Observation\IfNoValueF{#1}{: #1}}}
\DeclareDocumentCommand\addcitation{g}{\todo[size=\footnotesize,color = purple!10]{Add citation\IfNoValueF{#1}{: #1}}}

\theoremstyle{plain}
\newtheorem{theorem}{Theorem}[section]
\newtheorem*{theorem*}{Theorem}
\newtheorem{proposition}[theorem]{Proposition}
\newtheorem*{proposition*}{Proposition}

\newtheorem{corollary}[theorem]{Corollary}

\theoremstyle{definition}
\newtheorem{definition}[theorem]{Definition}

\theoremstyle{remark}

\newtheorem{remark}[theorem]{Remark}
\newtheorem{example}[theorem]{Example}
\newtheorem{notation}[theorem]{Notation}

\newcommand{\Mod}{\ensuremath{\mathsf{Mod}}}

\newcommand{\Obj}{\mathrm{Obj}}
\newcommand{\Lex}{\ensuremath{\mathsf{Lex}}}
\newcommand{\Rex}{\ensuremath{\mathsf{Rex}}}

\newcommand{\Spec}{\mathrm{Sp}}

\newcommand{\regcard}{\ensuremath{\mathsf{Reg}}}

\newcommand{\LC}{\ensuremath{\mathsf{LC}}}
\newcommand{\LP}{\ensuremath{\mathsf{LP}}}
\newcommand{\Groth}{\ensuremath{\mathsf{G}}}
\newcommand{\sites}{\ensuremath{\mathsf{S}}}
\newcommand{\op}{\ensuremath{\mathsf{op}}}
\newcommand{\co}{\ensuremath{\mathsf{co}}}

\newcommand{\id}{\ensuremath{1}}

\renewcommand{\lim}{\mathrm{lim}}
\newcommand{\colim}{\mathrm{colim}}

\newcommand{\AAA}{\mathfrak{a}}
\newcommand{\BBB}{\mathfrak{b}}

\newcommand{\CCC}{\mathfrak{c}}
\newcommand{\DDD}{\mathfrak{d}}
\newcommand{\EEE}{\mathfrak{e}}

\newcommand{\Set}{\ensuremath{\mathsf{Set}} }
\newcommand{\Ab}{\ensuremath{\mathsf{Ab}} }
\newcommand{\Cont}{\ensuremath{\mathsf{Cont}}}
\newcommand{\Cocont}{\ensuremath{\mathsf{Cocont}}}

\newcommand{\Sh}{\ensuremath{\mathsf{Sh}} }

\newcommand{\Qch}{\ensuremath{\mathsf{Qcoh}} }
\newcommand{\Coh}{\ensuremath{\mathsf{Coh}} }
\newcommand{\Cat}{\ensuremath{\mathsf{Cat}} }

\newcommand{\Psnat}{\ensuremath{\mathsf{Psnat}} }

\newcommand{\Ind}{\ensuremath{\mathsf{Ind}}}
\newcommand{\Fun}{\ensuremath{\mathsf{Fun}}}

\newcommand{\hlra}{\lhook\joinrel\longrightarrow}

\newcommand{\aaa}{\ensuremath{\mathcal{A}}}
\newcommand{\bbb}{\ensuremath{\mathcal{B}}}
\newcommand{\ccc}{\ensuremath{\mathcal{C}}}
\newcommand{\ddd}{\ensuremath{\mathcal{D}}}

\newcommand{\iii}{\ensuremath{\mathcal{I}}}
\newcommand{\jjj}{\ensuremath{\mathcal{J}}}

\newcommand{\LLL}{\ensuremath{\mathcal{L}}}

\newcommand{\ooo}{\ensuremath{\mathcal{O}}}

\newcommand{\uuu}{\ensuremath{\mathcal{U}}}
\newcommand{\vvv}{\ensuremath{\mathcal{V}}}

\newcommand{\calR}{\ensuremath{\mathscr{R}}}

\newcommand{\calT}{\ensuremath{\mathscr{T}}}

\title[Filtered bicolimit presentations and tensor products]{Filtered bicolimit presentations of locally presentable linear categories, Grothendieck categories and their tensor products}

\author{Julia Ramos Gonz\'alez}
\address[Julia Ramos Gonz\'alez]{Universit\'e Catholique de Louvain, Institut de Recherche en Math\'ematique et Physique, Chemin du Cyclotron 2/L7.01.02, 1348 Louvain-la-Neuve, Belgium.}
\email{\href{mailto:julia.ramos@uclouvain.be}{julia.ramos@uclouvain.be}}

\thanks{The author is a Postdoctoral Researcher of the Fonds de la Recherche
	Scientifique – FNRS [32709538]. This paper was partly written while she was a Postdoctoral Fellow of the Research Foundation – Flanders (FWO) [12T2619N]. She acknowledges as well the support of the Research Foundation Flanders (FWO) under Grant No. G.0112.13N during the time in which many of the results of this paper were obtained.}

\usepackage{pdfsync}
\begin{document}
	\begin{abstract}
		We investigate two different ways of recovering a Grothendieck category as a filtered bicolimit of small categories and the compatibility of both with the tensor product of Grothendieck categories. 
		Firstly, we show that any locally presentable linear category (and in particular any Grothendieck category) can be recovered as the filtered bicolimit of its subcategories of $\alpha$-presentable objects, with $\alpha$ varying in the family of small regular cardinals. We then prove that the tensor product of locally presentable linear categories (and in particular the tensor product of Grothendieck categories) can be recovered as a filtered bicolimit of the Kelly tensor product of $\alpha$-cocomplete linear categories of the corresponding subcategories of $\alpha$-presentable objects. Secondly, we show that one can recover any Grothendieck category as a filtered bicolimit of its linear site presentations. We then prove that the tensor product of Grothendieck categories, in contrast with the first case, cannot be recovered in general as a filtered bicolimit of the tensor product of the corresponding linear sites. 
		Finally, as a direct application of the first presentation, we translate the functoriality, associativity and symmetry of the Kelly tensor product of $\alpha$-cocomplete linear categories to the tensor product of locally presentable linear categories.
	\end{abstract}
\maketitle
\vspace{-1em}
\tableofcontents

\section{Introduction}\label{introduction} 

Grothendieck categories conform one of the most studied subfamilies within abelian categories, with the categories of modules and, more generally, the categories of quasi-coherent sheaves \cite[Tag 077P]{stacks-project} as the main and most well-known examples. They are essential objects in, among other fields, noncommutative algebraic geometry, where they play the role of (categorical models of) noncommutative schemes (see, for example \cite{noncommutative-projective-schemes}, \cite{noncommutative-curves-surfaces} or \cite{noncommutative-schemes} among many others). For convenience of the reader, we recall their definition:
\begin{definition}
	A \emph{Grothendieck category} is a cocomplete abelian category with a generator and exact filtered colimits.
\end{definition} 

The combination of the Gabriel-Popescu theorem \cite{caracterisation-categories-abeliennes-generateurs-limites-inductives-exactes} with the theory of enriched sheaves of Borceux and Quinteiro \cite{theory-enriched-sheaves} shows that Grothendieck categories are precisely the linear topoi, that is, Grothendieck categories are the categories of linear sheaves on small linear sites \cite{generalization-gabriel-popescu}. 

On the other hand, it is well-known that Grothendieck categories are in particular locally presentable \cite[Prop 3.4.16]{handbook-of-cat-alg3}. Therefore, thanks to the Gabriel-Ulmer duality \cite[Kor 7.11]{lokal-praesentierbare-kategorien}, we know that any Grothendieck category can be written as the $\Ind_{\alpha}$-completion of its full subcategory of $\alpha$-presentable objects (which is an $\alpha$-cocomplete linear category) for a big enough regular cardinal $\alpha$.

Both presentations, as categories of sheaves on a linear site or as $\Ind_{\alpha}$-completions, provide a powerful tool when working with Grothendieck categories. This is the case, for example, when defining the tensor product of Grothendieck categories, which is our main object of interest in this paper. 

The tensor product of Grothendieck categories $\boxtimes_{\Groth}$ \cite{tensor-product-linear-sites-grothendieck-categories} can be understood, in the framework of noncommutative algebraic geometry, as the correct notion of product between noncommutative schemes. Indeed, given two quasi-compact quasi-separated schemes $X,Y$ we have that 
\begin{equation}
	\Qch(X)\boxtimes_{\Groth} \Qch(Y) \cong \Qch(X \times Y)
\end{equation}
where $\Qch(-)$ denotes the category of quasi-coherent sheaves of a scheme and $\times$ is the (fibred) product of schemes (over $\Spec(\mathbb{Z})$). This equivalence was first proven for projective schemes in \cite[Thm 4.12]{tensor-product-linear-sites-grothendieck-categories}, and then, more generally, shown for quasi-compact quasi-separated schemes in \cite[Thm A]{localizations-tensor-categorires-fiber-products-schemes}.

The tensor product $\boxtimes_{\Groth}$ is defined in \cite{tensor-product-linear-sites-grothendieck-categories} in terms of the presentations of Grothendieck categories as categories of linear sheaves. More concretely, a tensor product of linear sites $\boxtimes_{\sites}$ is defined, and given two Grothendieck categories $\aaa, \bbb$ with sheaf representations $\aaa = \Sh(\AAA,\calT_{\AAA})$, $\bbb = \Sh(\BBB,\calT_\BBB)$ for linear sites $(\AAA,\calT_{\AAA}), (\BBB,\calT_\BBB)$, the tensor product $\aaa \boxtimes_{\Groth} \bbb$ is given by
\begin{equation}
	\aaa \boxtimes_{\Groth} \bbb = \Sh((\AAA,\calT_\AAA)\boxtimes_{\sites}(\BBB,\calT_\BBB)).
\end{equation}
This is well-defined, namely, the tensor product $\aaa \boxtimes_{\Groth} \bbb$ does not depend on the choice of linear sites providing the sheaf presentations of $\aaa$ and $\bbb$.

On the other hand, we also have a tensor product of locally presentable categories \cite{limits-2-categories-locally-presented-categories} at our disposal, which is easily translated to the $\Ab$-enriched setup, that is, to a tensor product of locally presentable linear categories. It is shown in \cite[Thm 5.4]{tensor-product-linear-sites-grothendieck-categories} that the tensor product of Grothendieck categories $\boxtimes_{\Groth}$  is an instance of the tensor product of locally presentable linear categories, which we denote by $\boxtimes_{\LP}$. 

The tensor product $\boxtimes_{\LP}$ is defined in \cite{limits-2-categories-locally-presented-categories} by means of the presentations of locally presentable categories as $\Ind_\alpha$-completions of $\alpha$-cocomplete small categories for suitable regular cardinal $\alpha$, and this construction works both for the nonenriched and the enriched setups. More concretely, in the linear case, we have that given $\aaa, \bbb$ two $\alpha$-locally presentable linear categories, the tensor product $\aaa \boxtimes_{\LP} \bbb$ is given by
\begin{equation}
	\aaa \boxtimes_{\LP} \bbb = \Ind_{\alpha}(\aaa^{\alpha} \otimes_{\alpha} \bbb^{\alpha})
\end{equation}
where $\aaa^{\alpha},\bbb^{\alpha}$ denote the corresponding linear subcategories of $\alpha$-presentable objects and $\aaa^{\alpha} \otimes_{\alpha} \bbb^{\alpha}$ is the Kelly-tensor product of $\alpha$-cocomplete linear categories  \cite{limits-2-categories-locally-presented-categories}, \cite{basic-concepts-enriched-category-theory}, \cite{structures-defined-finite-limits-enriched-context}.

The first goal of this paper is to show that linear site presentations and linear categories of $\alpha$-presentable objects provide us with two different ways to recover any Grothendieck category as a filtered bicolimit of small categories.

Our second goal is the understanding of the behaviour of these filtered bicolimits with respect to the tensor product of presentations. More concretely, given two locally presentable linear categories $\aaa,\bbb$, we show that the filtered bicolimit over the regular cardinals $\alpha$ of the Kelly tensor products of $\alpha$-cocomplete linear categories of the $\alpha$-presentable objects coincides with the tensor product $\aaa \boxtimes_{\mathsf{LP}} \bbb$. On the other hand, we show that the filtered bicolimit of the tensor product of linear site presentations does not allow us to recover the tensor product of the corresponding Grothendieck categories in general. 

Finally, as an application of the good behaviour of the filtered bicolimit presentations with respect to the Kelly tensor product of $\alpha$-cocomplete linear categories, we show that the functoriality, associativity and symmetry of the tensor product of Grothendieck categories can be also derived from the same properties of the Kelly tensor product. In particular, this allows us to make use of the functoriality of the Kelly tensor product in order to compute tensor products of cocontinuous linear functors between locally presentable linear categories and in particular, between Grothendieck categories.

Grothendieck categories are our main object of interest in this paper due to their relevance in noncommutative algebraic geometry. For this reason, we focus exclusively on linear enrichments, although the results we obtain most likely hold for other enrichments of Grothendieck topoi and locally presentable categories. In order to get a well-behaved theory of enriched Grothendieck topoi \cite{theory-enriched-sheaves} one needs to require that the category over which one enriches is regular in the sense of Barr, locally finitely presentable and symmetric monoidal closed, and moreover, that the full subcategory of finitely presentable objects is closed under the tensor product. We will not analyse here whether these assumptions on the enrichment are enough for our results to carry over.

\subsection*{Structure of the paper} 
We devote the first sections to the background and preliminary results required for our purposes:
\begin{itemize}
	\item[\textbf{In \S\ref{conventions-notations}:}] We fix the notation and conventions that will be used along the paper.
	\item[\textbf{In \S\ref{locallypresentable}:}] We revise the main aspects of the theory of locally presentable linear categories and their presentations as $\Ind_\alpha$-completions. In particular, we review the relation between the Kelly tensor product of $\alpha$-cocomplete linear categories and the tensor product of locally presentable linear categories and we study the behaviour of this relation when one raises the cardinality $\alpha$.
	\item[\textbf{In \S\ref{linearsites}:}] We revise the basic background on the theory of linear sites and their tensor product and how they present, respectively, Grothendieck categories and their tensor product.
	\item[\textbf{In \S\ref{parfilteredbicolim}:}] We revise briefly the general theory of $2$-filtered bicolimits from \cite{construction-2-filtered-bicolimits-categories}. In particular, we show that the instance of 2-filtered bicolimits in which we will be interested, namely, the filtered bicolimits (i.e. $2$-filtered bicolimits where the indexing category is just a filtered $1$-category), is well behaved with respect to the linear enrichement (see \Cref{klinear}). 
\end{itemize}

Once the necessary preliminaries are settled, we present our main results:
\begin{itemize}
	\item[\textbf{In \S\ref{pargrothasbicolim}:}] We provide two different ways of representing a Grothendieck category as a filtered bicolimit of small linear categories. The first one, generalizable to any locally presentable category, states the following:
	\begin{proposition}[{\Cref{thmgrothbicolim}}] \label{intropropbicolimit}
		Let $\ccc$ be a locally presentable linear category. Then $\ccc$ is a filtered bicolimit of its family of subcategories of locally $\alpha$-presentable objects $(\ccc^{\alpha})_{\alpha}$, where $\alpha$ varies in the totally ordered class of small regular cardinals. 
	\end{proposition}
	\noindent The second one makes use of the topos theoretical nature of Grothendieck categories. Given a Grothendieck category $\ccc$, we consider the category $\jjj_{\ccc}$ of all site presentations of $\ccc$, that is, all the LC morphisms $(\AAA,\calT) \to \ccc$ from a linear site $(\AAA,\calT)$ to $\ccc$ (see \Cref{defLC}), and we show it is filtered (see \Cref{jcfiltered}). We then consider the functor $G_{\ccc}: \jjj_{\ccc} \to \Cat$ assigning to each LC morphism its domain. We prove the following:
	\begin{theorem}[{\Cref{thmbicolimitsites}}]
	Given a Grothendieck category $\ccc$, we have that $\ccc$ is the linear filtered bicolimit of $G_{\ccc}$.
	\end{theorem}
	\item[\textbf{In \S\ref{tensorproduct}:}] We analyse the behaviour of the Kelly tensor product of $\alpha$-cocomplete linear categories and the tensor product of linear sites with respect to the corresponding filtered bicolimit presentations provided in \S\ref{pargrothasbicolim}. Given $\ccc,\ddd$ two locally $\kappa$-presentable linear categories, in virtue of \Cref{intropropbicolimit} above, one can show that the tensor product $\ccc \boxtimes_{\mathsf{LP}} \ddd$ of locally presentable categories can be recovered as the filtered bicolimit of the family $((\ccc \boxtimes_{\mathsf{LP}} \ddd)^{\alpha})_{\alpha \geq \kappa} = (\ccc^{\alpha} \otimes_{\alpha} \ddd^{\alpha})_{\alpha \geq \kappa}$ of subcategories of $\alpha$-presentable objects, with the transition functors given by the natural embeddings $(\ccc \boxtimes_{\mathsf{LP}} \ddd)^{\alpha} \subseteq (\ccc \boxtimes_{\LP} \ddd)^{\beta}$. However, at first sight, it is not clear whether these transition functors are precisely those induced, through the universal property of the Kelly tensor product, by the natural embeddings $\ccc^{\alpha} \subseteq \ccc^{\beta}$ and $\ddd^{\alpha} \subseteq \ddd^\beta$. We show that this is indeed the case:
	\begin{theorem}[{\Cref{productasbicolim}}]
		Let $\ccc,\ddd$ be two locally presentable linear categories. The tensor product $\ccc \boxtimes_{\mathsf{LP}} \ddd$ can be expressed as the filtered bicolimit of the tensor products $(\ccc^{\alpha} \otimes_{\alpha} \ddd^{\alpha})_{\alpha}$ of categories of $\alpha$-presentable objects, where $\alpha$ takes values in the totally ordered set of small regular cardinals, and the transition maps 
		$$\ccc^{\alpha} \otimes_{\alpha} \ddd^{\alpha} \to \ccc^{\beta} \otimes_{\beta} \ddd^{\beta}$$ 
		are those induced, through the universal property of the Kelly tensor product, by the canonical embeddings $\aaa^{\alpha} \subseteq \aaa^{\beta}$ and $\bbb^{\alpha} \subseteq \bbb^{\beta}$ for all $\alpha \leq \beta$.
	\end{theorem}
	\noindent While the Kelly tensor product of $\alpha$-cocomplete linear categories is well-behaved in this sense, this is not the case for the tensor product of linear sites. More concretely, given two Grothendieck categories $\aaa,\bbb$, we show in \Cref{remark:tplinearsites-vs-tpGrothendieck} that the filtered bicolimit of the tensor product of linear sites of the site presentations of $\aaa$ with the site presentations of $\bbb$ does not allow, in general, to recover $\aaa \boxtimes_{\Groth} \bbb$. 
	\item[\textbf{In \S\ref{parfunctassocsym}:}] We describe the functoriality of the tensor product of locally presentable linear categories via the functoriality of the Kelly tensor product of $\alpha$-cocomplete linear categories (see \Cref{functoriality}). This description provides an advantage when computing tensor products of cocontinuous linear functors in the cases where one has control of the subcategories of presentable objects. We illustrate this advantage by means of a geometrical example  (see \Cref{example:tp-qcoh}). To conclude, we describe the associativity and symmetry of the tensor product of locally presentable linear categories making use of the same properties of the Kelly tensor product of $\alpha$-cocomplete categories (see \Cref{propassociativity} and \Cref{propsymmetry}).
\end{itemize}  

\vspace{0,5cm}

\noindent \emph{Acknowledgements.} This article presents and extends part of the work carried out by the author in her PhD thesis under the supervision of Wendy Lowen and Boris Shoikhet. I am very grateful to both of them for the interesting discussions and their helpful comments. I would also like to thank an anonymous referee for the useful remarks that lead not only to the generalization of some results appearing in a previous version of the text but also to significant improvements on the presentation of the paper.

\section{Conventions and notation}\label{conventions-notations}
Throughout the paper we assume the axioms of the Tarski-Grothendieck (TG) set theory, that is, the axioms of the Zermelo-Fraenkel set theory with the axiom of choice (ZFC) together with the universe axiom. We fix a universe $\uuu$.

We denote by $\uuu$-$\regcard$ the totally ordered class of regular $\uuu$-small cardinals, which we will envision as a category. This category is not $\uuu$-small \cite[Cor 6.13]{Rathjen06}.

We fix a $\uuu$-small commutative ring $k$ for the rest of the article. We denote by $\uuu$-$\Cat(k)$ the $2$-category of $\uuu$-small $k$-linear categories (i.e. categories enriched over the category of $\uuu$-small $k$-modules $\uuu$-$\Mod(k)$) with $k$-linear functors and $k$-linear natural transformations. In particular, given $\AAA,\BBB \in \uuu$-$\Cat(k)$, we denote by $\Fun_k(\AAA,\BBB)$ the $\uuu$-small $k$-linear category of $k$-linear functors between $\AAA$ and $\BBB$. 

We denote by $\otimes = \otimes_k$ the \emph{tensor product of $k$-linear categories}, that is, given $\AAA,\BBB \in \uuu$-$\Cat(k)$, $\AAA \otimes \BBB$ is the $\uuu$-small $k$-linear category with objects the pairs $(A,B) \in \Obj(\AAA) \times \Obj(\BBB)$ and Hom-$k$-modules $\AAA\otimes \BBB((A,B),(A',B'))$ given by the tensor product of k-modules $\AAA(A,A') \otimes_k \BBB(B,B')$. The 2-category $\uuu$-$\Cat(k)$ endowed with $\otimes$ is a symmetric closed monoidal $2$-category and the inner hom corresponds with the external hom. More concretely, given $\AAA,\BBB,\CCC \in \uuu$-$\Cat(k)$, we have the universal property
	\begin{equation}\label{eq:tp-klinear-universalproperty}
		\Fun_k(\AAA \otimes_k \BBB, \CCC) \cong \Fun_k(\AAA, \Fun_k(\BBB,\CCC))
	\end{equation}  
	in $\uuu$-$\Cat(k)$.

Given $\AAA \in \uuu\mhyphen \Cat(k)$, we consider the \emph{category of $\uuu$-small (right) $\AAA$-modules} $\Fun_k(\AAA^{\op},\uuu$-$\Mod(k))$ and we denote it by $\uuu$-$\Mod(\AAA)$. Given $f \colon \AAA \to \BBB$ a morphism in $\uuu\mhyphen \Cat(k)$, we denote the \emph{restriction of scalars along $f$} by 
\begin{equation}
	f^*: \Mod(\BBB) \to \Mod(\AAA): F \mapsto F \circ f^\op, 
\end{equation}
its left adjoint, the \emph{extension of scalars along $f$}, by 
\begin{equation}
	f_{!}: \Mod(\AAA) \to \Mod(\BBB) 
\end{equation}
and its right adjoint, the \emph{coextension of scalars along $f$}, by 
\begin{equation}
	f_*: \Mod(\AAA) \to \Mod(\BBB).	
\end{equation}

Let $\alpha$ be a $\uuu$-small regular cardinal. We denote by $\uuu$-$\Cat_{\alpha}(k)$ the $2$-category of $\alpha$-cocomplete (i.e. closed under $\alpha$-small colimits) $\uuu$-small $k$-linear categories with $\alpha$-cocontinuous (i.e. $\alpha$-small colimit preserving) $k$-linear functors and $k$-linear natural transformations. In particular, given $\AAA,\BBB \in \uuu$-$\Cat_{\alpha}(k)$, we denote by $\Rex_{\alpha}(\AAA,\BBB)$ the $\uuu$-small $k$-linear category of $\alpha$-cocontinuous $k$-linear functors between $\AAA$ and $\BBB$.

For $\AAA \in \uuu$-$\Cat_{\alpha}(k)$, we consider the subcategory $\uuu$-$\Lex_{\alpha}(\AAA) \subseteq \uuu$-$\Mod(\AAA)$ of $\alpha$-continuous $\AAA$-modules, that is the $k$-linear functors $\AAA^{\op} \to \uuu$-$\Mod(k)$ preserving $\alpha$-small limits.

Given a small $k$-linear category $\AAA$, the categories $\uuu$-$\Lex_{\alpha}(\AAA)$ and $\uuu$-$\Mod(\AAA)$ are no longer $\uuu$-small, as neither is $\uuu$-$\regcard$, as pointed out above. However, one can always consider a larger universe $\vvv$ containing $\uuu$, such that $\uuu$-$\regcard$, $\uuu$-$\Lex_{\alpha}(\AAA)$ and $\uuu$-$\Mod(\AAA)$ are $\vvv$-small for all $\uuu$-small $k$-linear categories $\AAA$. From this point on and for the rest of the paper, we will omit the universes $\uuu$ and $\vvv$ from our notation and terminology.

To conclude, let us fix some notations we will use when working with $1$- and $2$-categories.

Given a $k$-linear $1$-category $\ccc$, we will denote the ($k$-module of) morphisms between two objects $A,B \in \ccc$ by $\ccc(A,B)$.

Let us in addition fix the notation we will use when working with 2-morphisms in the various bicategories (and 2-categories) appearing in this paper. In any bicategory, we denote by $\bullet$ the vertical composition of $2$-morphisms and by $\circ$ the horizontal composition of $2$-morphisms, following the convention in \cite{categories-working-mathematician}. In particular, given a diagram
\begin{equation}
	\begin{tikzcd}[column sep= 50pt ]
		A \arrow[r,bend left=40, "f"] \arrow[r,bend right=40,swap,"f"] \arrow[r, phantom,"\Downarrow \id_{f}", pos=0.6] &B  \arrow[r,bend left=40, "g"] \arrow[r,bend right=40, swap, "h"]  \arrow[r, phantom,"\Downarrow \alpha", pos=0.55] &C
	\end{tikzcd}
\end{equation}
in a bicategory $\ccc$, we denote by $\alpha \circ f$ to the horizontal composite $\alpha \circ \id_f$.

\section{Generalities on locally presentable linear categories}\label{locallypresentable}
In this section we provide an overview on the theory of locally presentable categories and their tensor product, focused on the aspects that will be needed along the text. We will work in the $k$-linear setup, which is almost identical to the non-enriched one, with the only essential difference being the fact that the monoidal structure in $\Set$ is cartesian, while the one in $\Mod(k)$ is not. We refer the reader to \cite{lokal-praesentierbare-kategorien} and \cite{locally-presentable-accessible-categories} for the non-enriched theory of locally presentable categories, and to \cite{structures-defined-finite-limits-enriched-context} for the enriched framework. 

%Recall that $\Cat(k)$ is a symmetric monoidal closed $2$-category. Given two small $k$-linear categories $\AAA, \BBB$, we denote by $\AAA \otimes \BBB$ its tensor product, which is the $k$-linear category with objects given by pairs of objects $(A,B) \in \AAA \times \BBB$ and morphisms $\AAA \otimes \BBB (A,B)$ given by the tensor product of morphisms $\AAA(A,A') \otimes_k \BBB(B,B')$. The internal hom is then given by the $k$-linear category $\Cat(k)(\AAA,\BBB)$ of $k$-linear functors with morphisms the $k$-linear natural transformations. We thus have
%\begin{equation}
%	\Cat(k)(\AAA \otimes \BBB, \CCC) \cong \Cat(k)(\AAA, \Cat(k)(\BBB,\CCC)).
%\end{equation} 

\subsection{Basic definitions and results}

\begin{definition}
	Let $\ccc$ be a $k$-linear category and $\alpha$ a small regular cardinal. We say that $C \in \ccc$ is $\alpha$-\emph{presentable} if the functor
	$$\ccc(C,-): \ccc \to \Mod(k)$$
	preserves $\alpha$-filtered colimits.  
\end{definition}

\begin{notation}
	Given a $k$-linear category $\ccc$, we denote by $\ccc^\alpha$ its full subcategory of $\alpha$-presentable objects. We denote by $\iota_{\alpha}: \ccc^\alpha \hlra \ccc$ the natural embedding. When necessary to avoid confusion, we will use a superscript to make explicit which category we are working with, e.g. $\iota^{\ccc}_{\alpha}: \ccc^\alpha \hlra \ccc$.
\end{notation}

\begin{proposition}[{\cite[Prop 1.16]{locally-presentable-accessible-categories}}]\label{proposition:alphasmallcolimits}
	Let $\ccc$ be a cocomplete category and $\alpha$ a small regular cardinal. Then $\ccc^\alpha$ is closed under $\alpha$-small colimits in $\ccc$. In particular, $\ccc^\alpha$ is $\alpha$-cocomplete and the embedding $\iota_{\alpha}: \ccc^\alpha \hlra \ccc$ is $\alpha$-cocontinuous.
\end{proposition}

\begin{remark}\label{remark:alphainbeta}
	Let $\ccc$ be $k$-linear category and $\alpha, \beta$ small regular cardinals. Assume that $\alpha < \beta$. Then, because $\beta$-filtered colimits are in particular $\alpha$-filtered, every $\alpha$-presentable object in $\ccc$ is also $\beta$-presentable and thus $\ccc^\alpha \subseteq \ccc^\beta$. Observe that, as a direct consequence of \Cref{proposition:alphasmallcolimits}, we have that $\ccc^\alpha$ is closed under $\alpha$-small colimits in $\ccc^\beta$. 
\end{remark}

\begin{notation}
	In the hypothesis of \Cref{remark:alphainbeta}, we denote by $\iota_{\alpha,\beta}: \ccc^{\alpha} \hlra \ccc^\beta$ the $\alpha$-cocontinuous natural embedding. When necessary to avoid confusion, we will use a superscript to make explicit which category we are working with, e.g. $\iota_{\alpha,\beta}^\ccc: \ccc^{\alpha} \hlra \ccc^\beta$.
\end{notation}	

\begin{definition}
	Let $\alpha$ be a small regular cardinal. A \emph{locally $\alpha$-presentable $k$-linear category} is a cocomplete $k$-linear category $\ccc$ with a small set of $\alpha$-presentable objects $\CCC$ such that every object in $\ccc$ is an $\alpha$-filtered colimit of elements in $\CCC$.
	
	We say that a category $\ccc$ is a \emph{locally presentable $k$-linear category} if there exists a small regular cardinal $\alpha$ such that $\ccc$ is a locally $\alpha$-presentable $k$-linear category.
\end{definition}

\begin{remark}[{\cite[Remark under Thm 1.20]{locally-presentable-accessible-categories}}]\label{remark:raising-cardinality}
	If a $k$-linear category $\ccc$ is locally $\beta$-presentable, then it is locally $\gamma$-presentable for all regular cardinals $\gamma \geq \beta$. Moreover, every object in $\ccc$ is $\alpha$-presentable for some regular cardinal $\alpha$, and thus we have that
	\begin{equation*}
		\ccc \cong \bigcup_{\alpha \in \regcard} \ccc^\alpha.
	\end{equation*}  
\end{remark}

\begin{notation}\,
	\begin{itemize}
		\item We denote by $\LP_\alpha(k)$ the $2$-category of locally $\alpha$-presentable $k$-linear categories with $1$-morphisms the cocontinuous $k$-linear functors preserving $\alpha$-presentable objects and $2$-morphisms the $k$-linear natural transformations. 
		\item We denote by $\LP(k)$ the $2$-category of locally presentable $k$-linear categories with $1$-morphisms the cocontinuous $k$-linear functors (i.e. the functors preserving all small colimits) and $2$-morphisms the $k$-linear natural transformations. In particular, given $\aaa,\bbb$ two locally presentable $k$-linear categories, we denote by $\Cocont(\aaa,\bbb)$ the $k$-linear category of cocontinuous $k$-linear functors between $\aaa$ and $\bbb$.
	\end{itemize}
\end{notation}

\begin{remark}[{\cite[Rem 1.19]{locally-presentable-accessible-categories}}]\label{rem:essentially-small}
	Given a locally $\alpha$-presentable $k$-linear category $\ccc$, its subcategory of $\alpha$-presentable objects $\ccc^{\alpha}$ is essentially small and hence we can consider it as an element in $\Cat_{\alpha}(k)$.
\end{remark}

Gabriel-Ulmer duality (see \cite[Kor 7.11]{lokal-praesentierbare-kategorien}) is a key result in the theory of locally presentable categories and we will use it repeatedly along the text. 

\begin{theorem}[Gabriel-Ulmer duality]\label{Gabriel-Ulmer}
	Let $\alpha$ be a small regular cardinal. There is a biequivalence of $2$-categories 
	\begin{equation*}
		(\Cat_{\alpha}(k))^{\co\op} \to \LP_\alpha(k),
	\end{equation*}
	where $(-)^{\co\op}$ denotes the conjugate-opposite $2$-category. 
\end{theorem}
Let's recall the explicit description of the duality. 
\begin{itemize}
	\item The biequivalence is given by the pseudofunctor 
	$$\Ind_\alpha: (\Cat_{\alpha}(k))^{\co\op} \to \LP_\alpha(k)$$
	sending an $\alpha$-cocomplete small category $\AAA$ to its $\Ind_{\alpha}$-completion $\Ind_{\alpha}(\AAA)$ (see \cite[\S3]{classification-accessible-categories}), which is locally $\alpha$-presentable \cite[Thm 1.46]{locally-presentable-accessible-categories}. In particular, we have that $\Ind_{\alpha}(\AAA) \cong \Lex_{\alpha}(\AAA)$ (see \cite[Thm 2.4]{classification-accessible-categories}) and the $k$-linear Yoneda embedding $Y_{\AAA}: \AAA \hookrightarrow \Mod(\AAA)$ factors as
	\begin{equation*}
		\begin{tikzcd}
			\AAA \arrow[r, "Y'_\alpha", hook] & \Lex_\alpha(\AAA) \arrow[r, "j_\alpha", hook] & \Mod(\AAA)
		\end{tikzcd}
	\end{equation*} 
	where $j_\alpha: \Lex_{\alpha}(\AAA) \hookrightarrow \Mod(\AAA)$ is the natural embedding. Moreover, the fully faithful functor $Y'_{\alpha}: \AAA \hookrightarrow \Lex_{\alpha}(\AAA)$ exhibits $\Lex_{\alpha}(\AAA)$ as the $\alpha$-free cocompletion of $\AAA$ \cite[Prop 1.45]{locally-presentable-accessible-categories}, that is, $Y'_{\alpha}$ preserves $\alpha$-small colimits and for every cocomplete linear category $\bbb$, we have that precomposition with $Y'_\alpha$ induces an equivalence
	\begin{equation}\label{eq:alpha-free-cocompletion}
		\Cocont(\Lex_{\alpha}(\AAA),\bbb) \cong \Rex_{\alpha}(\AAA,\bbb). 
	\end{equation}   
	\item The quasi-inverse pseudofunctor is given by
	$$(-)^\alpha: \LP_\alpha(k) \to (\Cat_{\alpha}(k))^{\co\op}$$
	which sends a locally $\alpha$-presentable category $\ccc$ to its subcategory of $\alpha$-presentable objects $\ccc^\alpha$, which is an $\alpha$-cocomplete (essentially) small category by \Cref{proposition:alphasmallcolimits} and \Cref{rem:essentially-small}.
\end{itemize}

\begin{remark}\label{remark:equivalenceLexalpha}
	Observe that, given a locally $\alpha$-presentable category $\ccc$, we have an equivalence of categories 
	$$E_{\alpha}: \ccc \overset{\cong}{\to} \Lex_{\alpha}(\ccc^\alpha).$$
	This equivalence is given on objects as follows:
	$$E_{\alpha}(C) \coloneqq \ccc(i_{\alpha}(-),C)$$
	where $\iota_{\alpha}: \ccc^{\alpha} \hookrightarrow \ccc$ is the natural embedding (see, e.g., the proof of \cite[Thm 1.46]{locally-presentable-accessible-categories}).
	%	Following \cite[Thm 1.46]{locally-presentable-accessible-categories}, we can describe the equivalence $E_\alpha$ explicitly as follows. Let $\iota_{\alpha}: \ccc^{\alpha} \hookrightarrow \ccc$ and $j_\alpha: \Lex_{\alpha}(\ccc^\alpha) \hookrightarrow \Mod(\ccc^\alpha)$ the natural embeddings. We have that the functor  
	%	$$Y_{\alpha}:\ccc \to \Mod(\ccc^{\alpha}): C \mapsto \ccc(i^{\ccc}_{\alpha}(-),C)$$ 
	%	factors through an equivalence
	%	\begin{equation*}
		%		\begin{tikzcd}
			%			\ccc \arrow[dr, "\cong", "E_{\alpha}^{\ccc}"'] \arrow[rr, "Y_{\alpha}"] &&\Mod(\ccc^{\alpha}).\\
			%			&\Lex_{\alpha}(\ccc^{\alpha}) \arrow[ur, hook, "j_\alpha"'] 
			%		\end{tikzcd}
		%	\end{equation*}
\end{remark}

\begin{proposition}[{\cite[Thm 1.46 \& Prop 1.35]{locally-presentable-accessible-categories}}]
	Let $\AAA \in \Cat_{\alpha}(k)$. Then, the natural embedding $j_\alpha: \Lex_{\alpha}(\AAA) \hookrightarrow \Mod(\AAA)$ preserves $\alpha$-filtered colimits and admits a right adjoint $R_\alpha: \Mod(\AAA) \to \Lex_\alpha(\AAA)$.
\end{proposition}

The following result will be useful for our purposes. 

\begin{proposition}\label{proposition:diagramlocallypresentable}
	Let $\ccc$ be a locally $\kappa$-presentable $k$-linear category and consider small regular cardinals $\beta > \alpha \geq \kappa$. Let $\iota_{\alpha,\beta}: \ccc^\alpha \hookrightarrow \ccc^\beta$ be the natural embedding functor. Then, the functor $(\iota_{\alpha,\beta})^*: \Mod(\ccc^\beta) \to \Mod(\ccc^\alpha)$ restricts to a functor
	$$(\iota_{\alpha,\beta})_s: \Lex_\beta(\ccc^\beta) \to \Lex_\alpha(\ccc^\alpha).$$
	Moreover, $(\iota_{\alpha,\beta})_s$ has a left adjoint $(\iota_{\alpha,\beta})^s: \Lex_{\alpha}(\ccc^\alpha) \to 
	\Lex_{\beta}(\ccc^\beta) $
	such that the diagram
	\begin{equation}\label{diagram:continuousLP}
		\begin{tikzcd}
			\ccc^\alpha \arrow[dd, "Y_{\ccc^\alpha}" description, hook] \arrow[rr, "{\iota_{\alpha,\beta}}" description, hook]                                &                                               & \ccc^\beta \arrow[dd, "Y_{\ccc^\beta}" description, hook]                                                                                           \\
			&                                               &                                                                                                                                                     \\
			\Mod(\ccc^\alpha) \arrow[dd, "R_\alpha" description] \arrow[rr, "{(\iota_{\alpha,\beta})_!}" description]          &                                               & \Mod(\ccc^\beta) \arrow[dd, "R_\beta" description] \\
			&                                               &                                                                                                                                                     \\
			\Lex_\alpha(\ccc^\alpha) \arrow[rr, "{(\iota_{\alpha,\beta})^s}" description, bend left] & \ccc \arrow[l, "\cong"', "E_\alpha"] \arrow[r, "\cong", "E_\beta"'] & \Lex_\beta(\ccc^\beta) 
		\end{tikzcd}
	\end{equation}
	is commutative up to isomorphism. In particular, the adjunction $(\iota_{\alpha,\beta})^s \dashv (\iota_{\alpha,\beta})_s$ is an equivalence of categories.
	\begin{proof}
		We first show that $(\iota_{\alpha,\beta})^*$ when restricted to the subcategory $\Lex_\beta(\ccc^\beta)$ takes values in $\Lex_{\alpha}(\ccc^\alpha)$. Indeed, if $X \in \Lex_{\beta}(\ccc^{\beta})$ one easily sees that
		\begin{equation*}
			(\iota_{\alpha,\beta})^*(X) = X((\iota_{\alpha,\beta})^\op(-))
		\end{equation*}
		preserves $\alpha$-small limits because both $(\iota_{\alpha,\beta})^\op$ and $X$ do. Therefore, $(\iota_{\alpha,\beta})^*$ does indeed restrict to a functor  
		$$(\iota_{\alpha,\beta})_s: \Lex_{\beta}(\ccc^\beta) \to \Lex_{\alpha}(\ccc^\alpha).$$ 
		More concretely, we have that the diagram 
		\begin{equation*}
			\begin{tikzcd}
				\Mod(\ccc^\alpha)                                                 &  & \Mod(\ccc^\beta) \arrow[ll, "{(\iota_{\alpha,\beta})^*}" description]                                               \\
				&  &                                                                                                                          \\
				\Lex_\alpha(\ccc^\alpha) \arrow[uu, "j_\alpha" description, hook] &  & \Lex_\beta(\ccc^\beta) \arrow[uu, "j_\beta" description, hook] \arrow[ll, "{(\iota_{\alpha,\beta})_s}" description]\\
				& \ccc \arrow[ul,"E_\alpha", "\cong"'] \arrow[ur,"E_\beta"', "\cong"]&
			\end{tikzcd}
		\end{equation*}
		is commutative, where the commutativity of the bottom triangle can be readily checked using the fact that $\iota_\beta \iota_{\alpha,\beta} = \iota_\alpha$. In particular, we have that $(\iota_{\alpha,\beta})_s$ is an equivalence.
		
		Define now the functor $(\iota_{\alpha,\beta})^s \coloneqq R_\beta (\iota_{\alpha,\beta})_! j_\alpha$. We show that this is the left adjoint to $(\iota_{\alpha,\beta})_s$. Indeed, we have that:
		\begin{equation}
			\begin{aligned}
				\Lex_\beta(\ccc^\beta) (R_\beta (\iota_{\alpha,\beta})_! j_\alpha (M), X) &\cong \Mod(\ccc^\beta)((\iota_{\alpha,\beta})_! j_\alpha (M), j_\beta(X))\\
				&\cong \Mod(\ccc^\alpha)(j_\alpha (M), (\iota_{\alpha,\beta})^*j_\beta(X)) \\
				&\cong \Mod(\ccc^\alpha)(j_\alpha (M), j_\alpha (\iota_{\alpha,\beta})_s(X))\\
				&\cong \Lex_\alpha(\ccc^\alpha)(M, (\iota_{\alpha,\beta})_s(X)).
			\end{aligned}
		\end{equation}
		Therefore, because $(\iota_{\alpha,\beta})_s$ was an equivalence, so is $(\iota_{\alpha,\beta})^s$, and it is in particular the quasi-inverse.
		
		It remains to check the commutativity of the diagram \eqref{diagram:continuousLP}. The commutativity up to isomorphism of the top square is a direct consequence of the ($k$-linear) Yoneda lemma. Indeed, given $C \in \ccc^\alpha$, we have that
		\begin{equation*}
			\begin{aligned}
				\Mod(\ccc^\beta)((\iota_{\alpha,\beta})_!Y_{\ccc^\alpha}(C), M) &\cong \Mod(\ccc^\alpha)(Y_{\ccc^\alpha}(C),(\iota_{\alpha,\beta})^*(M))\\
				&\cong (\iota_{\alpha,\beta})^*(M)(C)\\
				&= M(\iota_{\alpha,\beta}(C))\\
				&\cong \Mod(\ccc^\beta)(Y_{\ccc^\beta}(\iota_{\alpha,\beta}(C)),M)
			\end{aligned}
		\end{equation*} 
		for all $M \in \Mod(\ccc^\beta)$ and thus $(\iota_{\alpha,\beta})_!Y_{\ccc^\alpha}(C) \cong Y_{\ccc^\beta}\iota_{\alpha,\beta}(C)$ as desired. We now show the commutativity up to isomorphism of the triangle in the bottom. Using the fact that $(\iota_{\alpha,\beta})_s E_\beta = E_\alpha$, we obtain that
		\begin{equation*}
			E_\beta \cong (\iota_{\alpha,\beta})^s(\iota_{\alpha,\beta})_s E_\beta = (\iota_{\alpha,\beta})^sE_\alpha
		\end{equation*}
		as we wanted to prove. It remains to show the commutativity up to isomorphism of the bottom square. Because $(\iota_{\alpha,\beta})^s$ is an equivalence, it preserves small colimits, as so do $(\iota_{\alpha,\beta})_!, R_\alpha$ and $R_\beta$ because they are left adjoints. Therefore, as $\Mod(\ccc^\alpha)$ is generated via small colimits by the elements $\{Y_{\ccc^\alpha} (C)\}_{C \in \ccc^\alpha}$, we have that the bottom square in the diagram \eqref{diagram:continuousLP} is commutative up to isomorphism if and only if 
		\begin{equation*}
			R_\beta (\iota_{\alpha,\beta})_! Y_{\ccc^\alpha} \cong (\iota_{\alpha,\beta})^s R_\alpha Y_{\ccc^\alpha}.
		\end{equation*}
		Now observe that:
		\begin{equation*}
			\begin{aligned}
				(\iota_{\alpha,\beta})^s R_\alpha Y_{\ccc^\alpha} (C) &= R_\beta (\iota_{\alpha,\beta})_! j_\alpha R_\alpha Y_{\ccc^\alpha} (C)\\
				&\cong R_\beta (\iota_{\alpha,\beta})_! Y_{\ccc^\alpha} (C),
			\end{aligned}	
		\end{equation*}
		where the second isomorphism comes from the fact that $Y_{\ccc^\alpha} = j_\alpha Y'_{\ccc^\alpha}$ and $R_\alpha j_\alpha \cong \id_{\Lex_\alpha(\ccc^\alpha)}$. This concludes the proof. 
	\end{proof}
\end{proposition}

\begin{remark}
	\Cref{proposition:diagramlocallypresentable} can be compared to the topos-theoretic result \cite[Exp.III, Prop 1.2]{SGA4-1} (see \cite[Prop 2.3]{grothendieck-categories-bilocalization-linear-sites} for the $k$-linear counterpart). In fact, our choice of the notation $(-)^s$ and $(-)_s$ is borrowed from \cite{SGA4-1}. 
\end{remark}

\subsection{The tensor product: $\alpha$-cocomplete linear categories and locally presentable linear categories}
In this section we revise the monoidal structures in the $2$-categories $\Cat_{\alpha}(k)$, $\LP_{\alpha}(k)$ and $\LP(k)$ and the relations among them.

The $2$-category $\Cat_{\alpha}(k)$ of $\alpha$-cocomplete $k$-linear categories can be endowed with a closed monoidal structure \cite[\S 6.5]{basic-concepts-enriched-category-theory} as follows:
\begin{definition}\label{definition:alpha-tensor-product}
	Given $\AAA, \BBB \in \Cat_{\alpha}(k)$, there exists an $\AAA \otimes_{\alpha} \BBB \in \Cat_{\alpha}(k)$ and a $k$-linear functor $u^{\alpha}_{\AAA,\BBB}: \AAA \otimes \BBB \rightarrow \AAA \otimes_{\alpha} \BBB$ which is $\alpha$-cocontinuous in each variable, such that, for every $\CCC \in \Cat_{\alpha}(k)$, composition with $u^\alpha_{\AAA,\BBB}$ induces an equivalence
	\begin{equation}\label{alphaprod}
		\Rex_\alpha(\AAA \otimes_{\alpha} \BBB, \CCC) \cong \Rex_{\alpha,\alpha}(\AAA \otimes \BBB, \CCC) \cong \Rex_\alpha(\AAA, \Rex_\alpha(\BBB, \CCC)) 
	\end{equation} 
	in $\Cat_{\alpha}(k)$, where $\Rex_{\alpha,\alpha}(\AAA \otimes \BBB, \CCC)$ denotes the full subcategory of $\Fun_k(\AAA \otimes \BBB, \CCC)$ consisting of the functors that are $\alpha$-cocontinuous in each variable. We say that $\AAA \otimes_{\alpha} \BBB$ is the \emph{$\alpha$-Kelly tensor product} of $\AAA$ and $\BBB$ \cite{basic-concepts-enriched-category-theory},\cite{structures-defined-finite-limits-enriched-context}.
\end{definition} 
\begin{proposition}
	The $2$-category $\Cat_{\alpha}(k)$ together with the $\alpha$-Kelly tensor product $\otimes_\alpha$ is a closed monoidal bicategory. The monoidal unit is given by the full subcategory $\left( \Mod(k) \right)^\alpha \subseteq \Mod(k)$ of $\alpha$-presentable $k$-modules and the inner hom between $\AAA$ and $\BBB$ is given by $\Rex_{\alpha}(\AAA,\BBB)$.
\end{proposition}

\begin{remark}\label{remark:alphaprod}
	Given $\AAA, \BBB$ two $\alpha$-cocomplete $k$-linear categories and $\ccc$ a cocomplete $k$-linear category, one has that
	\begin{equation}\label{eq:alphaprod}
		\Rex_\alpha(\AAA \otimes_{\alpha} \BBB, \ccc) \cong \Rex_\alpha(\AAA, \Rex_\alpha(\BBB, \ccc)) \cong \Rex_{\alpha,\alpha}(\AAA \otimes \BBB, \ccc).
	\end{equation} 
\end{remark}

The $\alpha$-Kelly tensor product can be explicitly constructed as follows:
\begin{proposition}[{\cite[Prop 6.21]{basic-concepts-enriched-category-theory}}]\label{prop:explicit-construction-Kelytp}
	Consider $\AAA$, $\BBB \in \Cat_{\alpha}(k)$. The category $\AAA \otimes_{\alpha} \BBB$ can be recovered as the closure under $\alpha$-small colimits of the image of the composite 
	\begin{equation}\label{eq:universal-functor-alpha-kellytp}
		\begin{tikzcd}
			\AAA \otimes \BBB \arrow[r, "Y_{\AAA \otimes \BBB}", hook] & \Mod(\AAA \otimes \BBB) \arrow[r, "{R_{\alpha,\alpha}}"] & {\Lex_{\alpha,\alpha}(\AAA,\BBB)}
		\end{tikzcd}
	\end{equation}
	where $\Lex_{\alpha,\alpha}(\AAA,\BBB)$ denotes the full subcategory of $\Mod(\aaa^\alpha \otimes \bbb^\alpha)$ consisting of the bimodules $F:(\aaa^\alpha)^{\op}\otimes (\bbb^\alpha)^{\op} \to \Mod(k)$ that are $\alpha$-continuous in each variable, and the functor $R_{\alpha,\alpha}: \Mod(\AAA \otimes \BBB) \to \Lex_{\alpha,\alpha}(\AAA,\BBB)$ is the left adjoint to the natural embedding $j_{\alpha,\alpha}:\Lex_{\alpha,\alpha}(\AAA,\BBB) \hookrightarrow \Mod(\AAA \otimes \BBB)$. 
	
	Moreover, the universal functor $u^\alpha_{\AAA,\BBB}: \AAA \otimes \BBB \to \AAA \otimes_{\alpha} \BBB$ is given by the corestriction of the functor \eqref{eq:universal-functor-alpha-kellytp}.
\end{proposition}

The closed monoidal structure on $\Cat_{\alpha}(k)$ can be transported via the Gabriel-Ulmer duality (\Cref{Gabriel-Ulmer}) to a closed monoidal structure in $\LP_{\alpha}(k)$ (see \cite[\S5.1]{limits-2-categories-locally-presented-categories}). More concretely, we have the following statement:
\begin{proposition}\label{prop:tp-alpha-presentable}
	The 2-category $\LP_{\alpha}(k)$ is a closed monoidal bicategory with the bifunctor $\boxtimes_{\LP_\alpha}:\LP_{\alpha}(k) \otimes \LP_{\alpha}(k) \to \LP_{\alpha}(k)$ given by
	\begin{equation*}
		\aaa \boxtimes_{\LP_\alpha} \bbb \coloneqq \Lex_\alpha(\aaa^\alpha \otimes_\alpha \bbb^\alpha).
	\end{equation*} 
	The monoidal unit is given by $\Mod(k)$ and the inner hom between $\aaa$ and $\bbb$ is given by $\Ind_{\alpha}(\Rex_\alpha(\aaa^\alpha,\bbb^\alpha))$. 
\end{proposition}

\begin{remark}\label{rem:presentable-objects-tensor-product}
	Given $\aaa, \bbb$ locally $\alpha$-presentable categories, we have that
	\begin{equation*}
		(\aaa \boxtimes_{\LP_\alpha} \bbb)^\alpha \cong \aaa^\alpha \otimes_\alpha \bbb^\alpha. 
	\end{equation*}	
\end{remark}

\begin{remark}\label{rem:alpha-in-each-variable}
	Observe that, by substituting $\ccc$ in \Cref{eq:alphaprod} with $\Mod(k)^\op$, one obtains that 
	\begin{equation*}
		\aaa \boxtimes_{\LP_\alpha} \bbb \cong \Lex_{\alpha,\alpha}(\aaa^\alpha,\bbb^\alpha).
	\end{equation*}	
\end{remark}

In \cite[\S 5.1]{limits-2-categories-locally-presented-categories} it is shown that the tensor product $\boxtimes_{\LP_\alpha}$ admits a nice description independent of the cardinal $\alpha$.
\begin{proposition}\label{proposition:Bird}
	The tensor product $\boxtimes_{\LP_\alpha}:\LP_{\alpha}(k) \otimes \LP_{\alpha}(k) \to \LP_{\alpha}(k)$ is isomorphic to the bifunctor
	\begin{equation*}
		\LP_{\alpha}(k) \otimes \LP_{\alpha}(k) \to \LP_{\alpha}(k): (\aaa , \bbb)  \mapsto \Cont(\aaa^\op,\bbb).
	\end{equation*}
\end{proposition} 

By virtue of \Cref{remark:raising-cardinality} and \Cref{proposition:Bird}, one can extend the tensor product in $\LP_{\alpha}(k)$ to a tensor product in $\LP(k)$ \cite[\S 5.1]{limits-2-categories-locally-presented-categories}.
\begin{proposition}\label{prop:tplocallypresentable}
	The bifunctor $\boxtimes_{\LP}:\LP(k) \otimes \LP(k) \to \LP(k)$ given by
	\begin{equation*}
		\aaa \boxtimes_{\LP} \bbb \coloneqq \Cont (\aaa^\op,\bbb)
	\end{equation*} 
	induces a closed monoidal structure in $\LP(k)$, with monoidal unit given by $\Mod(k)$ and inner hom given by $\Cocont(\bbb,\ccc)$. More concretely, given $\aaa,\bbb \in \LP(k)$, there exists a functor $s_{\aaa,\bbb}: \aaa \otimes \bbb \to \aaa \boxtimes_{\LP} \bbb$ which is cocontinuous in each variable and such that precomposition with $s^\alpha_{\aaa,\bbb}$ provides an equivalence
	\begin{equation}\label{eq:universal-property-tplp}
		\Cocont(\aaa \boxtimes_{\LP} \bbb,\ccc) \cong \Cocont,\Cocont(\aaa \otimes \bbb, \ccc) \cong \Cocont(\aaa \otimes \bbb, \ccc) 
	\end{equation} 
	for all $\ccc \in \LP(k)$, where $\Cocont,\Cocont(- \otimes -, -)$ denotes the $k$-linear functors that are cocontinuous in each variable. 
	Moreover, if both $\aaa$ and $\bbb$ are $\alpha$-presentable one has that
	\begin{equation*}
		(\aaa \boxtimes_{\LP} \bbb)^\alpha \cong \aaa^\alpha \otimes_\alpha \bbb^\alpha. 
	\end{equation*}	
\end{proposition}

\begin{remark}\label{rem:compatibility-universal-properties}
	The universal properties of $\otimes_\alpha$ and $\boxtimes_\LP$ are compatible in a nice way (see the proof of \cite[Prop 5.3]{limits-2-categories-locally-presented-categories}). In other words, if $\aaa$ and $\bbb$ are locally $\alpha$-presentable $k$-linear categories, the diagram
	\begin{equation*}
		\begin{tikzcd}
			\aaa^\alpha \otimes \bbb^\alpha \arrow[d, "{u_{\aaa^\alpha,\bbb^\alpha}}"'] \arrow[rrr, "i^\aaa_\alpha \otimes i^\bbb_\alpha" description, hook] &  &  & \aaa \otimes \bbb \arrow[d] \arrow[dd, "{s_{\aaa,\bbb}}", bend left=80]            \\
			\aaa^\alpha \otimes_\alpha \bbb^\alpha \arrow[rrr, "Y'_{\aaa^\alpha \otimes_\alpha \bbb^\alpha}" description, hook] \arrow[d, "\cong" description]          &  &  & \Lex_\alpha(\aaa^\alpha \otimes_\alpha \bbb^\alpha) \arrow[d, "\cong" description] \\
			(\aaa \boxtimes_\LP \bbb)^\alpha \arrow[rrr, "i^{\aaa \boxtimes_\LP \bbb}_\alpha" description, hook]                                                        &  &  & \aaa \boxtimes_\LP \bbb                                                           
		\end{tikzcd}
	\end{equation*}
	is commutative.
\end{remark}

\begin{notation}\label{not:alpha-alpha-vs-alpha}
	Let $\aaa$ and $\bbb$ be locally $\alpha$-presentable categories. We know by virtue of \Cref{prop:tp-alpha-presentable} and \Cref{rem:alpha-in-each-variable} that $\aaa \boxtimes_{\LP} \bbb \cong \Lex_{\alpha,\alpha}(\aaa^\alpha,\bbb^\alpha)$ is an $\alpha$-presentable category and that its subcategory of $\alpha$-presented objects is given by $\aaa^\alpha \otimes_{\alpha} \bbb^\alpha$. We denote by $k_{\alpha}$ the natural embedding $\aaa^\alpha \otimes_{\alpha} \bbb^\alpha \hookrightarrow \Lex_{\alpha,\alpha}(\aaa^\alpha,\bbb^\alpha)$ (see \Cref{prop:explicit-construction-Kelytp}). Therefore, we have that the functor $$\Lex_{\alpha,\alpha}(\aaa^\alpha,\bbb^\alpha) \to \Mod(\aaa^\alpha \otimes_{\alpha} \bbb^\alpha): X \mapsto \Lex_{\alpha,\alpha}(\aaa^\alpha,\bbb^\alpha) (k_{\alpha}(-), X)$$ 
	factors via an equivalence as follows
	\begin{equation*}
		\begin{tikzcd}
			{\Lex_{\alpha,\alpha}(\aaa^\alpha,\bbb^\alpha)} \arrow[rr] \arrow[rd, "\cong" description, "E_\alpha"'] &                                                                      & \Mod(\aaa^\alpha \otimes_\alpha \bbb^\alpha) \\
			& \Lex_\alpha(\aaa^\alpha \otimes_\alpha \bbb^\alpha) \arrow[ru, hook, "j_\alpha"'] &                                             
		\end{tikzcd}
	\end{equation*}
\end{notation}
Moreover, we have that the following diagram is commutative by definition
\begin{equation*}
	\begin{tikzcd}[row sep=large]
		\aaa^\alpha \otimes \bbb^\alpha \arrow[d, "R_{\alpha,\alpha} Y_{\aaa^\alpha \otimes \bbb^\alpha}" description] \arrow[rr, "{u_{\aaa^\alpha,\bbb^\alpha}}" description] &  & \aaa^\alpha \otimes_{\alpha} \bbb^\alpha \arrow[lld, "k_\alpha" description, hook] \arrow[d, "Y'_{\aaa^\alpha \otimes \bbb^\alpha}" description, hook] \\
		{\Lex_{\alpha,\alpha}(\aaa^\alpha,\bbb^\alpha)} \arrow[rr, "E_\alpha" description]                                                                     &  & \Lex_\alpha(\aaa^\alpha \otimes_{\alpha} \bbb^\alpha)                                                       
	\end{tikzcd}
\end{equation*}

%In addition, we know that given locally $\alpha$-presentable categories $\aaa, \bbb$, we have that $\aaa \boxtimes_{\mathsf{LP}} \bbb = \Lex_{\alpha}(\aaa^{\alpha},\bbb^{\alpha})$ is $\alpha$-locally presentable and its subcategory of $\alpha$-presentable objects is given by $\aaa^{\alpha} \otimes_{\alpha} \bbb^{\alpha}$. For these results, we point the reader to \cite{structures-defined-finite-limits-enriched-context} and \cite{basic-concepts-enriched-category-theory}, or to \cite{tensor-product-finitely-cocomplete-abelian-categories} for the case $\alpha=\aleph_0$. 

\subsubsection{Raising the cardinality}
Let $\aaa$, $\bbb$ be two locally presentable $k$-linear categories and choose the smallest regular cardinal $\kappa$ such that both are locally $\kappa$-presentable. Consider regular cardinals $\alpha \leq \beta$ and denote by $\iota^{\aaa}_{\alpha,\beta}: \aaa^{\alpha} \hookrightarrow \aaa^{\beta}$, $\iota^{\bbb}_{\alpha,\beta}: \bbb^{\alpha} \hookrightarrow \bbb^{\beta}$ the natural embeddings, which in particular are $\alpha$-cocontinuous. Observe that we have a canonical $\alpha$-cocontinuous morphism
\begin{equation}\label{transitionfunctors}
	f_{\alpha, \beta}: \aaa^{\alpha} \otimes_\alpha \bbb^{\alpha} \rightarrow \aaa^{\beta} \otimes_\beta \bbb^{\beta}
\end{equation}
that makes the diagram 
\begin{equation}\label{eq:raising-cardinality-Kelly}
	\begin{tikzcd}
		\aaa^{\alpha} \otimes \bbb^{\alpha} \arrow[r,"\iota^{\aaa}_{\alpha,\beta} \otimes \iota^{\bbb}_{\alpha,\beta}"] \arrow[d,"u_{\aaa^{\alpha},\bbb^{\alpha}}"'] &\aaa^{\beta} \otimes \bbb^{\beta} \arrow[d,"u_{\aaa^{\beta},\bbb^{\beta}}"]\\
		\aaa^{\alpha} \otimes_{\alpha} \bbb^{\alpha} \arrow[r,"f_{\alpha,\beta}"] &\aaa^{\beta} \otimes_{\beta} \bbb^{\beta}.
	\end{tikzcd}
\end{equation}
commutative.
Indeed, $f_{\alpha,\beta}$ is defined as the image via the universal property (\ref{alphaprod}) in $\Rex_\alpha(\aaa^{\alpha} \otimes_{\alpha} \bbb^{\alpha}, \aaa^{\beta} \otimes_{\beta} \bbb^{\beta})$ of the composite $$\aaa^{\alpha} \otimes \bbb^{\alpha} \xrightarrow{\iota^{\aaa}_{\alpha,\beta} \otimes \iota^{\bbb}_{\alpha,\beta}} \aaa^{\beta} \otimes \bbb^{\beta} \xrightarrow{u_{\aaa^{\beta}, \bbb^{\beta}}} \aaa^{\beta} \otimes_{\beta} \bbb^{\beta},$$ which is $\alpha$-cocontinuous in each variable.

The following remark will be useful for the proof of \Cref{proposition:2variablesdiagramlocallypresentable} below.
\begin{remark}
	Let $\alpha, \beta$ be two regular cardinals and $\aaa, \bbb$ two $\gamma$-presentable categories where $\gamma = \min(\alpha,\beta)$. Then the tensor-inner hom adjunction of $k$-linear categories from \eqref{eq:tp-klinear-universalproperty} restricts to an isomorphism
	\begin{equation} \label{eq:innerhomalphabeta}
		\Phi_{\alpha,\beta}: \Lex_{\alpha,\beta}(\aaa^\alpha,\bbb^\beta) \to \Lex_{\alpha}((\aaa^\alpha)^\op, \Lex_\beta(\bbb^\beta)).
	\end{equation} 
	More concretely, a functor $F \in \Lex_{\alpha,\beta}(\aaa^\alpha,\bbb^\beta)$ is sent to the functor $\Phi_{\alpha,\beta}(F) \in \Lex_{\alpha}((\aaa^\alpha)^\op, \Lex_\beta(\bbb^\beta))$ defined by $\Phi_{\alpha,\beta}(F)(A)(B) \coloneqq F(A,B)$ for all $A \in \aaa^\alpha, B \in \bbb^\beta$ and a natural transformation $\sigma: F \Rightarrow G$ in $\Lex_{\alpha,\beta}(\aaa^\alpha,\bbb^\beta)$ is sent to the natural transformation $\Phi_{\alpha,\beta}(\sigma): \Phi_{\alpha,\beta}(F) \Rightarrow \Phi_{\alpha,\beta}(G)$ defined by $\left[ \left[ \Phi_{\alpha,\alpha}(\sigma)\right]_A\right]_B \coloneqq \sigma_{(A,B)}$ for all $A \in \aaa^\alpha, B \in \bbb^\beta$. We will denote by
	\begin{equation}\label{eq:innerhomalphabetainverse}
		\Psi_{\alpha,\beta}: \Lex_{\alpha}((\aaa^\alpha)^\op, \Lex_\beta(\bbb^\beta))  \to \Lex_{\alpha,\beta}(\aaa^\alpha,\bbb^\beta)
	\end{equation}
	the inverse of this isomorphism. More concretely, any $F \in \Lex_{\alpha}((\aaa^\alpha)^\op, \Lex_{\beta}(\bbb^\beta))$ is sent to the functor $\Psi_{\alpha,\beta}(F) \in \Lex_{\alpha,\beta}(\aaa^\alpha,\bbb^\beta)$ given by $\Psi_{\alpha,\beta}(F)(A,B) = F(A)(B)$ and a natural transformation $\eta: F \Rightarrow G$ is sent to the natural transformation $\Psi_{\alpha,\beta}(\eta): \Psi_{\alpha,\beta}(F) \Rightarrow \Psi_{\alpha,\beta}(G)$ defined by $\left[\Psi_{\alpha,\beta}(\eta) \right]_{(A,B)} = \left[ \eta_A\right]_B$.
\end{remark}

The following is a $2$-variables version of \Cref{proposition:diagramlocallypresentable}.
\begin{proposition}\label{proposition:2variablesdiagramlocallypresentable}
	Let $\aaa, \bbb$ be locally $\kappa$-presentable $k$-linear categories and consider small regular cardinals $\beta > \alpha \geq \kappa$. Let $\iota \coloneqq \iota^\aaa_{\alpha, \beta} \otimes \iota^\bbb_{\alpha, \beta} : \aaa^\alpha \otimes \bbb^\alpha \to  \aaa^\beta \otimes \bbb^\beta$ be the natural functor. Then, the functor $\iota^*: \Mod(\aaa^\beta \otimes \bbb^\beta) \to \Mod(\aaa^\alpha \otimes \bbb^\alpha)$ sends objects in $\Lex_{\beta,\beta}(\aaa^\beta,\bbb^\beta)$ to objects in $\Lex_{\alpha,\alpha}(\aaa^\alpha, \bbb^\alpha)$ and hence restricts to a functor
	$$\iota_s: \Lex_{\beta,\beta}(\aaa^\beta,\bbb^\beta) \to \Lex_{\alpha,\alpha}(\aaa^\alpha, \bbb^\alpha).$$
	Moreover, $\iota_s$ has a left adjoint $\iota^s: \Lex_{\alpha,\alpha}(\aaa^\alpha, \bbb^\alpha) \to 
	\Lex_{\beta,\beta}(\aaa^\beta,\bbb^\beta)$
	such that the diagram
	\begin{equation}\label{diagram:2variables-continuousLP}
		\begin{tikzcd}
			\aaa^\alpha \otimes \bbb^\alpha \arrow[rr, "\iota" description] \arrow[dd, hook]                                                              &  & \aaa^\beta \otimes \bbb^\beta \arrow[dd, hook]                                                                                       \\
			&  &                                                                                                                                      \\
			\Mod(\aaa^\alpha \otimes \bbb^\alpha) \arrow[dd, "R_{\alpha,\alpha}" description] \arrow[rr, "\iota_!" description]      &  & \Mod(\aaa^\beta \otimes \bbb^\beta) \arrow[dd, "R_{\beta,\beta}" description]    \\
			&  &                                                                                                                                      \\
			{\Lex_{\alpha,\alpha}(\aaa^\alpha,\bbb^\alpha)} \arrow[rr, "\iota^s" description]  &  & {\Lex_{\beta,\beta}(\aaa^\beta,\bbb^\beta)} 
		\end{tikzcd}
	\end{equation}
	is commutative up to isomorphism. Furthermore, the adjunction $\iota^s \dashv \iota_s$ is an equivalence of categories.
	\begin{proof}
		We first prove that $\iota^s$ when restricted to the subcategory $\Lex_{\beta,\beta}(\aaa^\beta, \bbb^\beta)$ takes values in $\Lex_{\alpha,\alpha}(\aaa^\alpha, \bbb^\alpha)$. Indeed, given $X \in \Lex_{\beta,\beta}(\aaa^\beta, \bbb^\beta)$, we have that 
		\begin{equation*}
			\begin{aligned}
				\iota^*(X) (-,-) = X (\iota^\op(-,-)) = X \left( (\iota^\aaa_{\alpha, \beta})^\op (-), (\iota^\bbb_{\alpha, \beta})^\op(-)\right) 
			\end{aligned}
		\end{equation*}
		belongs to $\Lex_{\alpha,\alpha}(\aaa^\alpha, \bbb^\alpha)$ because $(\iota^\aaa_{\alpha, \beta})^\op$ and $(\iota^\bbb_{\alpha, \beta})^\op$ preserve $\alpha$-small limits and $X$ preserves $\beta$-small (and in particular $\alpha$-small) limits in each variable. We thus have a functor $\iota_s:\Lex_{\beta,\beta}(\aaa^\beta,\bbb^\beta) \to \Lex_{\alpha,\alpha}(\aaa^\alpha,\bbb^\alpha)$ such that the diagram
		\begin{equation*}
			\begin{tikzcd}
				\Mod(\aaa^\alpha \otimes \bbb^\alpha)                                                               &  & \Mod(\aaa^\beta \otimes \bbb^\beta) \arrow[ll, "\iota^*" description]                                                           \\
				&  &                                                                                                                                 \\
				{\Lex_{\alpha,\alpha}(\aaa^\alpha,\bbb^\alpha)} \arrow[uu, "{j_{\alpha,\alpha}}" description, hook] &  & {\Lex_{\beta,\beta}(\aaa^\beta,\bbb^\beta)} \arrow[uu, "{j_{\beta,\beta}}" description, hook] \arrow[ll, "\iota_s" description]
			\end{tikzcd}
		\end{equation*}
		is commutative.
		Define now the functor $\iota^s \coloneqq R_{\beta,\beta} \iota_! j_{\alpha,\alpha}$. We are going to show that $\iota^s$ is a left adjoint of $\iota_s$. Indeed, given $M \in \Lex_{\alpha,\alpha}(\aaa^\alpha,\bbb^\alpha)$ and $N \in \Lex_{\beta,\beta}(\aaa^\beta,\bbb^\beta)$ we have that
		\begin{equation*}
			\begin{aligned}
				\Lex_{\beta,\beta}(\aaa^\beta,\bbb^\beta)\left(\iota^s(M), N \right) &\cong \Mod(\aaa^\beta \otimes \bbb^\beta)\left(\iota_! j_{\alpha,\alpha} (M), j_{\beta,\beta} (N) \right)\\
				&\cong \Mod(\aaa^\alpha \otimes \bbb^\alpha)\left(j_{\alpha,\alpha} (M), \iota^* j_{\beta,\beta} (N) \right)\\
				&\cong \Mod(\aaa^\alpha \otimes \bbb^\alpha)\left(j_{\alpha,\alpha} (M), j_{\alpha,\alpha} \iota_s (N) \right)\\
				&\cong \Lex_{\alpha,\alpha}(\aaa^\alpha,\bbb^\alpha)\left(M, \iota_s(N) \right),
			\end{aligned}
		\end{equation*}
		naturally in $M$ and $N$.
		We now prove that the diagram \eqref{diagram:2variables-continuousLP} is commutative up to isomorphism. The commutativity of the top square follows from the Yoneda lemma by the same argument as the one used in the proof of \Cref{proposition:diagramlocallypresentable}. We now prove that the bottom square of the diagram commutes up to isomorphism.
		Observe that, for all $X \in \Mod(\aaa^\alpha \otimes \bbb^\alpha)$ and $M \in \Lex_{\beta,\beta}(\aaa^\beta,\bbb^\beta)$ we have that
		\begin{equation*}
			\begin{aligned}
				\Lex_{\beta,\beta}(\aaa^\beta,\bbb^\beta)\left(\iota^s R_{\alpha,\alpha}(X), M \right) &\cong \Lex_{\alpha,\alpha}(\aaa^\alpha,\bbb^\alpha)\left(R_{\alpha,\alpha}(X), \iota_s(M) \right)\\
				&\cong \Mod(\aaa^\alpha \otimes \bbb^\alpha)\left( X, j_{\alpha,\alpha} \iota_s (M)\right)\\
				&\cong \Mod(\aaa^\alpha \otimes \bbb^\alpha)\left( X, \iota^* j_{\beta,\beta} (M)\right)\\  
				&\cong \Mod(\aaa^\beta \otimes \bbb^\beta)\left(\iota_! (X), j_{\beta,\beta}(M) \right)\\
				&\cong \Lex_{\beta,\beta}(\aaa^\beta,\bbb^\beta)\left(R_{\beta,\beta} \iota_! (X), M \right),
			\end{aligned}
		\end{equation*}
		naturally in $X$ and $M$ and thus by Yoneda lemma we can conclude that $R_{\beta,\beta} \iota_! \cong \iota^s R_{\alpha,\alpha}$ as desired.
		
		It remains to show that the adjunction $\iota^s \dashv \iota_s$ is an equivalence adjunction. A direct inspection shows that, $\iota_s = - \circ \iota^\op = - \circ ((\iota_{\alpha,\beta}^\aaa)^\op \otimes (\iota_{\alpha,\beta}^\bbb)^\op)$ coincides with the following composite of equivalences:
		\settowidth{\LabelWidth}{\WidestLabel}
		\begin{equation*}
			\begin{aligned}
				\Lex_{\beta,\beta}(\aaa^\beta,\bbb^\beta) &\LabelCong{\Phi_{\beta,\beta}}  \Lex_{\beta}((\aaa^\beta)^\op, \Lex_\beta(\bbb^\beta))\\
				&\LabelCong{(\iota_{\alpha,\beta}^\bbb)_s \circ -} \Lex_{\beta}((\aaa^\beta)^\op, \Lex_\alpha(\bbb^\alpha))\\
				&\LabelCong{\Psi_{\beta,\alpha}} \Lex_{\beta,\alpha}(\aaa^\beta,\bbb^\alpha)\\
				&\LabelCong{S_{\aaa^\beta}^{\bbb^\alpha}} \Lex_{\alpha,\beta}(\bbb^\alpha, \aaa^\beta)\\
				&\LabelCong{\Phi_{\alpha,\beta}} \Lex_{\alpha}((\bbb^\alpha)^\op, \Lex_\beta(\aaa^\beta))\\
				&\LabelCong{(\iota_{\alpha,\beta}^\aaa)_s \circ -} \Lex_{\alpha}((\bbb^\alpha)^\op, \Lex_\alpha(\aaa^\alpha))\\
				&\LabelCong{\Psi_{\alpha,\alpha}} \Lex_{\alpha,\alpha}(\bbb^\alpha,\aaa^\alpha)\\
				&\LabelCong{S_{\bbb^\alpha}^{\aaa^\alpha}} \Lex_{\alpha,\alpha}(\aaa^\alpha,\bbb^\alpha),
			\end{aligned}
		\end{equation*} 
		where $\Phi_{-,-}$ and $\Psi_{-,-}$ are the isomorphisms defined in \eqref{eq:innerhomalphabeta} and \eqref{eq:innerhomalphabetainverse} respectively, $(\iota_{\alpha,\beta}^\bbb)_s$ and $(\iota_{\alpha,\beta}^\aaa)_s$ are the equivalences obtained from applying \Cref{proposition:diagramlocallypresentable} to $\bbb$ and $\aaa$ respectively, and $S_{\aaa^\beta}^{\bbb^\alpha}$ and $S_{\bbb^\alpha}^{\aaa^\alpha}$ are given by the symmetry of the tensor product in $\Cat(k)$, i.e. by the braidings $S_{\aaa^\beta}^{\bbb^\alpha}: \aaa^\beta \otimes \bbb^\alpha \cong \bbb^\alpha \otimes \aaa^\beta$ and $S_{\bbb^\alpha}^{\aaa^\alpha}:\bbb^\alpha \otimes \aaa^\alpha \cong \aaa^\alpha \otimes \bbb^\alpha$, respectively. Therefore $\iota_s$ is an equivalence, and hence so is $\iota^s$, as we wanted to prove.
	\end{proof}
\end{proposition}

\section{Generalities on linear sites}\label{linearsites}

In this section we revise the basic notions and results concerning $k$-linear sites, as they will be an essential tool in the rest of the paper. For a more complete account we point the reader to \cite[\S 2]{linearized-topologies-deformation-theory} and \cite[\S2]{grothendieck-categories-bilocalization-linear-sites}.	

\subsection{Basic definitions and results}
$k$-Linear sites and Grothendieck $k$-linear categories can be seen as the $k$-linear counterpart of the classical Grothendieck sites and Grothendieck topoi from \cite{SGA4-1}. We point the reader to \cite{theory-enriched-sheaves} for more general enhancements of sites and topoi.

Let $\AAA$ be a small $k$-linear category.

\begin{definition}
	Given an object $A \in \AAA$, a \emph{$k$-linear sieve on $A$} is a subobject $R$ of the representable module $\AAA(-,A)$ in the category $\Mod(\AAA)$. Given a family $S = (f_i: A_i \rightarrow A)_{i \in I}$ of morphisms in $\AAA$, the \emph{$k$-linear sieve generated by $S$} is the smallest $k$-linear sieve $R$ on $A$ such that $f_i \in R(A_i)$ for all $i \in I$. 
\end{definition}

\begin{definition}
	A \emph{$k$-linear cover system} $\mathscr{R}$ on $\AAA$ consists of providing for each $A \in \AAA$ a family of $k$-linear sieves $\calR(A)$ on $A$. The $k$-linear sieves in a $k$-linear cover system $\calR$ are called \emph{covering sieves} or \emph{covers} (for $\calR$). We will say that a family $(f_i: A_i \to A)_{i\in I}$ is a \emph{cover}, or a \emph{covering family}, if the $k$-linear sieve it generates is a cover.
\end{definition}

\begin{definition}
	Given $R$ a $k$-linear sieve on $A \in \AAA$ and $g : A' \rightarrow A$ a morphism in $\AAA$, the \emph{pullback of $R$ along $g$}, denoted $g^{-1}R$, is the $k$-linear sieve on $A'$ obtained as the pullback
	\begin{equation*}
		\begin{tikzcd}
			g^{-1}R \arrow[r, hook] \arrow[d] &\AAA(-,A') \arrow[d,"g \circ \,-"] \\
			R \arrow[r,hook] &\AAA(-,A)
		\end{tikzcd}
	\end{equation*}
	in $\Mod(\AAA)$. In particular, we have that $g^{-1}R (A'') = \{ f : A'' \rightarrow A' \,\,|\,\, g \circ f \in R(A'')\}$ for all $A'' \in \AAA$.
\end{definition}

\begin{definition}
	A $k$-linear cover system $\calT$ on $\AAA$ is \emph{localizing} if it satisfies the following:
	\begin{itemize}	[leftmargin=40pt]
		\item[\textbf{(Id)}] Identity axiom: Given any object $A \in \AAA$, the $k$-linear sieve generated by $\id_A$ is a cover for $\calT$, i.e. $\AAA(-,A) \in \calT(A)$ for all $A \in \AAA$; 
		\item[\textbf{(Pb)}] Pullback axiom: Given a covering sieve $R \subseteq \calT(A)$ and $g : A' \rightarrow A$ a morphism in $\AAA$, the pullback sieve $g^{-1}R$ is also a covering sieve for $\calT$. 
	\end{itemize}
	If moreover $\calT$ also satisfies the following:
	\begin{itemize}[leftmargin=40pt]
		\item[\textbf{(Glue)}] Glueing axiom: Let $R$ be a $k$-linear sieve on $A$. If there exists a $k$-linear sieve $S$ on $A$ such that for all morphism $g: A' \rightarrow A$ in $S$ the pullback sieve $g^{-1}R \in \calT(A')$, then $R \in \calT(A)$;
	\end{itemize}
	we say $\calT$ is a \emph{$k$-linear Grothendieck topology}.
\end{definition}	

\begin{definition}
	A \emph{$k$-linear site} is a pair $(\AAA,\calT)$ where $\AAA$ is a $k$-linear category and $\calT$ is a $k$-linear Gro\-then\-dieck topology on $\AAA$.
\end{definition}

Given a $k$-linear site $(\AAA,\calT)$ one can define linearised versions of presheaves and sheaves, in analogy with the classical notions.

\begin{definition}
	\begin{itemize}
		\item[ ] 
		\item A \emph{$k$-linear presheaf} $F$ on $(\AAA, \calT)$ is an $\AAA$-module, this is $F \in \Mod(\AAA)$.
		\item A \emph{$k$-linear sheaf} $F$ on $(\AAA, \calT)$ is a $k$-linear presheaf such that the restriction functor
		\begin{equation*}
			F(A) \cong \Mod(\AAA)(\AAA(-,A), F) \to \Mod(\AAA)(R, F)
		\end{equation*}
		is an isomorphism for all $A \in \AAA$ and all $R \in \calT(A)$. We denote by $\Sh(\AAA, \calT) \subseteq \Mod(\AAA)$ the full subcategory of $k$-linear sheaves.
	\end{itemize}	
\end{definition}

\begin{definition}
	Consider a small $k$-linear category $\AAA$. A $k$-linear Grothendieck topology $\calT$ on $\AAA$ is called \emph{subcanonical} if for every $A \in \AAA$, the representable presheaf $\AAA(-,A) \in \Mod(\AAA)$ is a sheaf, that is, it belongs to $\Sh(\AAA,\calT)$. The largest Grothendieck subcanonical topology is called the \emph{canonical topology}.
\end{definition}

\begin{remark}\label{canonicaltop}
	We will frequently consider Grothendieck $k$-linear categories $\ccc$ themselves as (large) $k$-linear sites, endowed with their canonical topology $\calT_{\ccc,\mathrm{can}}$. In this particular case, the covering families are the jointly epimorphic families and $\Sh(\ccc, \calT_{\ccc,\mathrm{can}}) \cong \ccc$.
\end{remark}

The following is a consequence of Gabriel-Popescu theorem \cite{caracterisation-categories-abeliennes-generateurs-limites-inductives-exactes} in combination with enriched topos theory \cite{theory-enriched-sheaves}.
\begin{theorem}
	The categories of $k$-linear sheaves over $k$-linear sites are precisely the Grothendieck $k$-linear categories.
\end{theorem}

%The fact that Grothendieck $k$-linear categories are in particular locally presentable $k$-linear categories \cite[Prop 3.4.16]{handbook-of-cat-alg3}, allows us to combine the (generalization of) Gabriel-Popescu theorem as presented \cite[Thm 3.7]{generalization-gabriel-popescu} with the representation theorem of locally presentable enriched categories \cite[Thm 7.2 + \S7.4]{structures-defined-finite-limits-enriched-context} to obtain the following.
%
%\begin{theorem}\label{factorizationgrothendieck}
%Let $\ccc$ be an $\alpha$-presentable Grothendieck $k$-linear category. Then the functor
%	$$Y_{\alpha}:\ccc \to \Mod(\ccc^{\alpha}): C \mapsto \ccc(i^{\ccc}_{\alpha}(-),C),$$ 
%	where $\ccc^{\alpha}$ denotes the full subcategory of $\ccc$ of $\alpha$-presentable objects, factors through an equivalence
%	\begin{equation*}
	%	\begin{tikzcd}
		%	\ccc \arrow[dr, "\cong"', "E_{\alpha}^{\ccc}"] \arrow[rr, "Y_{\alpha}"] &&\Mod(\ccc^{\alpha}).\\
		%	&\Lex_{\alpha}(\ccc^{\alpha}) = \Sh(\ccc^{\alpha}) \arrow[ur, hook] 
		%	\end{tikzcd}
	%	\end{equation*}
%\end{theorem}

\begin{definition}\label{defcontinuous}
	Let $(\AAA,\calT_{\AAA})$ and $(\BBB,\calT_{\BBB})$ be $k$-linear sites. We say that a $k$-linear functor $f:\AAA \to \BBB$ is a \emph{continuous morphism of $k$-linear sites} if the restriction of scalars functor $f^*:\Mod(\BBB) \to \Mod(\AAA)$ preserves $k$-linear sheaves. 
	
	We denote by $f_s: \Sh(\BBB,\calT_{\BBB}) \to \Sh(\AAA, \calT_{\AAA})$ the corresponding restriction functor, and by $f^s: \Sh(\AAA, \calT_{\AAA}) \to \Sh(\BBB,\calT_{\BBB})$ its left adjoint.
\end{definition}

The class of LC morphisms between sites, where LC stands for ``Lemme de comparison'' \cite[\S 4]{generalization-gabriel-popescu}, will be used in the following sections: 

\begin{definition}[{\cite[Def 3.4]{tensor-product-linear-sites-grothendieck-categories}}]\label{defLC}
	Consider a $k$-linear functor $f: \AAA \to \CCC$, where $\AAA$ is small and $\CCC$ may be large. 
	\begin{enumerate}
		\item Suppose $\CCC$ is endowed with a $k$-linear cover system $\calT_{\CCC}$. We say that $f: \AAA \to (\CCC, \calT_{\CCC})$ satisfies 
		\begin{itemize}
			\item[\textbf{(G)}] if for every $C \in \CCC$ there is a covering family $(f(A_i) \to C)_{i \in I}$ for $\calT_{\CCC}$.
		\end{itemize}
		\item Suppose $\AAA$ is endowed with a $k$-linear cover system $\calT_{\AAA}$. We say that $f: (\AAA, \calT_{\AAA}) \to \CCC$ satisfies
		\begin{itemize}
			\item[\textbf{(F)}] if for every morphism $c: f(A) \to f(A')$ in $\CCC$ there exists a covering family $(a_i: A_i \to A)_{i \in I}$ for $\calT_{\AAA}$ and morphisms $f_i: A_i \to A'$ for each $i \in I$ such that $cf(a_i) = f({f}_i)$ for all $i \in I$;
			\item[\textbf{(FF)}] if for every morphism $a: A \to A'$ in $\AAA$ with $f(a) = 0$ there exists a covering family $(a_i: A_i \to A)_{i \in I}$ for $\calT_{\AAA}$ with $aa_i = 0$ for all $i \in I$. 
		\end{itemize}
		\item Suppose $\AAA$ and $\CCC$ are endowed with $k$-linear cover systems $\calT_{\AAA}$ and $\calT_{\CCC}$ respectively. We say that $f: (\AAA, \calT_{\AAA}) \to (\CCC, \calT_{\CCC})$ satisfies
		\begin{itemize}
			\item[\textbf{(LC)}] if $f$ satisfies (G) with respect to $\calT_{\CCC}$, (F) and (FF) with respect to $\calT_{\AAA}$, and we further have $\calT_{\AAA} = f^{-1} \calT_{\CCC}$.
		\end{itemize}
	\end{enumerate}
	We denote the class of LC morphisms by $\LC$.
\end{definition}

\begin{remark}\label{propclassLC}
	Observe that the identity morphism of a $k$-linear site $(\AAA,\calT)$ belongs to $\LC$. In addition, $\LC$ is closed under composition \cite[Prop 4.4]{generalization-gabriel-popescu}.
\end{remark}

The importance of LC morphisms between linear sites relies in the fact that they are continuous morphisms inducing equivalences between the corresponding sheaf categories \cite[Cor 4.5]{generalization-gabriel-popescu} together with the following key result:
\begin{theorem}[{\cite[Thm 5]{grothendieck-categories-bilocalization-linear-sites}}]\label{roofthm}
	Let $(\AAA, \calT_{\AAA})$ and $(\BBB, \calT_{\BBB})$ be $k$-linear sites and consider a colimit preserving functor $F: \Sh(\AAA,\calT_{\AAA}) \to \Sh(\BBB,\calT_{\BBB})$. There exist a subcanonical site $(\CCC,\calT_{\CCC})$ and a diagram
	\begin{equation*}
		\begin{tikzcd}
			&\CCC\\
			\AAA \arrow[ur,"f"] &&\BBB, \arrow[ul,"w",swap]
		\end{tikzcd}
	\end{equation*}
	where $f$ is a continuous morphism of $k$-linear sites and $w$ is an LC morphism, such that
	\begin{equation*}
		\begin{tikzcd}
			\Sh(\AAA, \calT_{\AAA}) \arrow[rr,"F"] \arrow[dr,"f^s",swap] &&\Sh(\BBB, \calT_{\BBB})\\
			&\Sh(\CCC,\calT_{\CCC}) \arrow[ur,"w_s",swap]
		\end{tikzcd}
	\end{equation*}
	is a commutative diagram up to isomorphism. 
\end{theorem}	

\subsection{The tensor product: linear sites and Grothendieck categories}	
In \cite{tensor-product-linear-sites-grothendieck-categories} a tensor product $\boxtimes_{\sites}$  of Grothendieck categories was introduced based on the definition of a tensor product of $k$-linear sites. 
\begin{definition}
	Given two $k$-linear sites $(\AAA ,\calT_{\AAA}), (\BBB,\calT_{\BBB})$, the \emph{tensor product of $k$-linear sites} $(\AAA ,\calT_{\AAA}) \boxtimes_{\sites} (\BBB,\calT_{\BBB})$ is given by the $k$-linear site $(\AAA \otimes \BBB, \calT_\AAA \boxtimes_{\sites} \calT_\BBB)$, where $\calT_\AAA \boxtimes_{\sites} \calT_\BBB$ is the smallest Grothendieck topology on $\AAA \otimes \BBB$ such that, for each object $(A,B) \in \AAA \otimes \BBB$, $\calT_\AAA \boxtimes_{\sites} \calT_\BBB(A,B)$ contains the images of all the canonical morphisms
		\begin{equation*}
			R \otimes_k S \to \AAA(-,A) \otimes_k \BBB(-,B) = \AAA \otimes \BBB (-,(A,B))
		\end{equation*}
	where $R \in \calT_{\AAA}(A)$ and $S \in \calT_{\BBB}(B)$.
\end{definition}

The tensor product of Grothendieck $k$-linear categories is defined in \cite{tensor-product-linear-sites-grothendieck-categories} as follows.

\begin{definition}
	Given two Grothendieck $k$-linear categories $\aaa = \Sh(\AAA,\calT_{\AAA})$ and $\bbb = \Sh(\BBB,\calT_{\BBB})$, the \emph{tensor product of Grothendieck $k$-linear categories} $\aaa \boxtimes_{\Groth} \bbb$  is given by
	$$\aaa \boxtimes_{\Groth} \bbb = \Sh(\AAA \otimes \BBB, \calT_{\AAA} \boxtimes_{\sites} \calT_{\BBB})$$
\end{definition}

\begin{remark}
	The tensor product is well-defined, namely, it is independent of the $k$-linear site presentations of $\aaa$ and $\bbb$ chosen \cite[Prop 4.1]{tensor-product-linear-sites-grothendieck-categories}. The proof of this relies on the following result regarding LC morphisms.
\end{remark}

\begin{proposition}[{\cite[Prop 3.14]{tensor-product-linear-sites-grothendieck-categories}}]\label{tpLC}
	Let $f: (\AAA,\calT_{\AAA}) \to (\BBB,\calT_{\BBB})$ and $g: (\CCC,\calT_{\CCC}) \to (\DDD,\calT_{\DDD})$ be $k$-linear functors. If $f$ and $g$ are LC morphisms, so is the functor 
	$$f \otimes g: (\AAA \otimes \CCC,\calT_{\AAA} \boxtimes_{\sites} \calT_{\CCC}) \to (\BBB \otimes \DDD,\calT_{\BBB} \boxtimes_{\sites} \calT_{\DDD}).$$
\end{proposition}

The following result provides a useful description of the tensor product of Grothendieck $k$-linear categories.

\begin{proposition}[{\cite[Cor 2.16 \& Cor 2.18]{tensor-product-linear-sites-grothendieck-categories}}]
	Consider Grothendieck $k$-linear categories $\aaa = \Sh(\AAA,\calT_{\AAA})$ and $\bbb = \Sh(\BBB,\calT_{\BBB})$. Then, the tensor product $\aaa \boxtimes_{\Groth} \bbb = \Sh(\AAA \otimes \BBB,\calT_{\AAA} \boxtimes_{\sites} \calT_{\BBB}) \subseteq \Mod(\AAA \otimes \BBB)$ is given by the full subcategory of bimodules $F:\AAA^{\op} \otimes \BBB^{\op} \to \Mod(k)$ for which $F(A,-) \in \Sh(\BBB,\calT_{\BBB})$ for all $A \in \AAA$ and $F(-,B) \in \Sh(\AAA,\calT_{\AAA})$ for all $B \in \BBB$.
\end{proposition} 

Recall that Grothendieck categories are in particular locally presentable categories \cite[Prop 3.4.16]{handbook-of-cat-alg3}. The following result was proven in \cite[Thm 5.4]{tensor-product-linear-sites-grothendieck-categories}.
\begin{theorem}
	Given two Grothendieck categories $\aaa, \bbb$, we have that 
	$$\aaa \boxtimes_{\Groth} \bbb \cong \aaa \boxtimes_{\mathsf{LP}} \bbb.$$
\end{theorem}

\section{Generalities on the 2-filtered bicolimit of categories}\label{parfilteredbicolim}
Given a (pseudo)functor $F: \aaa \rightarrow \Cat$ where $\aaa$ is a category and $\Cat$ denotes the $2$-category of small categories, we can consider \emph{Grothendieck's construction} of the colimit category $\underrightarrow{\lim}_{\aaa}(F)$ \cite[Expos\'e VI]{SGA4-2}. In particular, this construction can be performed when $\aaa$ is a filtered category. In \cite{construction-2-filtered-bicolimits-categories} a suitable generalization of the Grothendieck construction to the 2-categorical realm is provided for the filtered case, and referred to as \emph{$2$-filtered bicolimit}. In this section we provide a short overview of the $2$-filtered bicolimit of loc.cit. More concretely, we focus on its construction for the particular 1-categorical case of a pseudofunctor $F: \aaa \to \Cat$ where $\aaa$ is a filtered $1$-category. In addition, we show that in such case, if the functor $F$ factors through the $2$-category of $k$-linear categories $\Cat(k)$, the filtered bicolimit is also $k$-linear. 

We recall some important definitions.

\begin{definition}
	Let $F,G: \aaa \to \bbb$ be two pseudofunctors between $2$-categories $\aaa, \bbb$. 
	\begin{enumerate}
		\item A \emph{pseudonatural transformation} $\Phi: F \Rightarrow G$ is given by a family of $1$-morphisms $$(\Phi_A: F(A) \to G(A))_{A \in \Obj(\aaa)}$$ and a family of invertible $2$-morphisms $$(\Phi_f: \Phi(A') \circ F(f) \Rightarrow G(f) \circ \Phi(A))_{(f:A \to A') \in \aaa}$$ with the corresponding coherence laws, which we do not write down here. 
		\item Given two pseudonatural transformations $\Phi,\Psi: F \Rightarrow G$, a \emph{morphism of pseudonatural transformations} $r$ between $\Phi$ and $\Psi$ is a modification, that is a family of $2$-morphisms $$(r_A: \Phi_A \Rightarrow \Psi_A)_{A \in \aaa}$$ such that $$(G(f)\circ r_A)\bullet \Phi_f = \Psi_f \bullet (r_{A'} \circ F(f))$$ for all $f:A \to A'$ in $\aaa$.
	\end{enumerate}
	We denote by $\Psnat(F,G)$ the category of pseudonatural transformations between $F$ and $G$, with morphisms given by the modifications (see \cite{construction-2-filtered-bicolimits-categories} or \cite{basic-bicategories}).
\end{definition}

The notion of \emph{2-filtered 2-category} is introduced in \cite[\S2]{construction-2-filtered-bicolimits-categories} as a suitable generalisation in the $2$-categorical realm of the classical notion of filtered category. In particular, as it is already mentioned in the introduction of \cite{construction-2-filtered-bicolimits-categories}, any ($1$-)category considered as a trivial 2-category is 2-filtered if and only if it is filtered as an ordinary category. Throughout this paper we will always use an indexing category which is of this latter type, hence we can safely avoid going through the technicalities of the construction of $2$-filtered bicolimit for a general indexing $2$-category. 

%We are now in position to define the $2$-filtered bicolimit \cite[Thm. 1.19]{construction-2-filtered-bicolimits-categories}.
\begin{definition}[{\cite[Thm. 1.19]{construction-2-filtered-bicolimits-categories}}]\label{defbicolimit}
	Given a pseudofunctor $F: \iii \rightarrow \Cat$, where $\iii$ is a $2$-filtered $2$-category\footnote{The bicolimit and its construction actually work when the indexing category is a \emph{pre $2$-filtered $2$-category}, which is a weaker notion than that of $2$-filtered $2$-category, as pointed out in the introduction of \cite{construction-2-filtered-bicolimits-categories}.}, the \emph{2-filtered bicolimit of $F$} is a category $\LLL(F)$ together with a pseudonatural transformation $F \Rightarrow \LLL(F)$ from $F$ to the constant $2$-functor $\iii \rightarrow \Cat$ taking the value $\LLL(F)$ such that, for every category $\ccc$, it induces via composition an equivalence of categories
	\begin{equation}\label{univprop}
	\Cat(\LLL(F), \ccc) \cong \Psnat(F, \ccc)
	\end{equation}
	between the category of functors $\LLL(F) \rightarrow \ccc$ and the category of pseudonatural transformations between the pseudofunctor $F$ and the constant $2$-functor taking the value $\ccc$.
\end{definition}
\begin{remark}
	Note that such a category is uniquely determined up to a unique equivalence.
\end{remark}
\begin{remark}
	Observe that the original definition (see \cite[Thm. 1.19]{construction-2-filtered-bicolimits-categories}) only considers 2-filtered bicolimits of $F$ with $F$ a strict $2$-functor. For our purposes we need to consider the more general situation in which $F$ is a pseudofunctor.
\end{remark}

The main result of \cite{construction-2-filtered-bicolimits-categories} provides, given a 2-functor $F: \iii \to \Cat$, a construction of the bicolimit $\LLL(F)$ in an intrinsic way in terms of the $2$-functor $F$, which in particular makes of the equivalence \eqref{univprop} an isomorphism of categories. One can observe that, when the indexing category is just an ordinary filtered ($1$-)category the construction is greatly simplified. 
%In particular, we will be focusing on ordinary filtered categories given by a directed poset. 
For this choice of indexing category, one can easily extend the construction from \cite{construction-2-filtered-bicolimits-categories} to the case in which $F$ is a pseudofunctor by means of a slight generalization of the results explained in loc.cit.\footnote{The general construction for $F$ a pseudofunctor should be possible in full generality, without restrictions of the indexing category, just by readjusting the notion of \emph{homotopy} in \cite[1.5(iii)]{construction-2-filtered-bicolimits-categories}, as we have done in our particular case. } 

We flesh out below the construction of the bicolimit for our particular situation, i.e. when
\begin{itemize}
	\item $F: \iii \to \Cat$ is not necessarily a strict 2-functor but a pseudofunctor,
	\item the indexing category $\iii$ is a filtered (1-)category.
\end{itemize}
 
%Consider $\iii$ an ordinary filtered category induced by a directed poset (in particular, between two objects there is at most one morphism). Consider a pseudofunctor $F: \iii \rightarrow \Cat$. We define the following category $\LLL(F)$:

\noindent\textbf{Description of the objects:}

Objects of $\LLL(F)$ are pairs $(x,A)$ where $A \in \iii$ and $x \in F(A)$.

\noindent\textbf{Description of the morphisms:}

First we describe the class of premorphisms:
\begin{itemize}
	\item  Consider two objects $(x, A)$ and $(y, B)$. A \emph{premorphism} $(x,A) \rightarrow (y,B)$ consists of a triple $(u,f,v)$, where $u: A \rightarrow C$, $v: B \rightarrow C$ and $f: F(u)(x) \rightarrow F(v)(y)$ in $F(C)$. In order to make the object $C$ explicit in the notation, we will write $(u,f,v)_C: (x,A) \to (y,B)$. 
	\item Two pre-morphisms $(u_1,f,v_1)_{C_1}, (u_2,g,v_2)_{C_2}:(x,A) \rightarrow (y,B)$ are said to be \emph{homotopic} if there exists an object $C \in \iii$ and morphisms $w_i: C_i \rightarrow C$, for $i =1,2$ such that $w_1 \circ u_1 = w_2 \circ u_2$, $w_1 \circ v_1 = w_2 \circ v_2$ and the following diagram commutes
	\begin{equation*} \rule{\leftmargin}{0in}
		\begin{tikzcd}[column sep=0.8em]
		F(w_1)\circ F(u_1)(x) \arrow[r, "\cong"] \arrow[d, "F(w_1)(f_1)"'] & F(w_1 \circ u_1)(x) = F(w_2 \circ u_2)(x) \arrow[r, "\cong"] & F(w_2) \circ F(u_2)(x) \arrow[d, "F(w_2)(f_2)"] \\
		F(w_1) \circ F(v_1)(y) \arrow[r, "\cong"] & F(w_1 \circ v_1)(y) = F(w_2 \circ v_2)(y) \arrow[r, "\cong"] & F(w_2)\circ F(v_2)(y).
		\end{tikzcd}
	\end{equation*}
	This relation is an equivalence relation and we denote the equivalence class of a premorphism $(u,f,v)_C:(x,A) \to (y,B)$ by $[(u,f,v)_C]$. 
\end{itemize}
By means of the homotopy relation, we define the morphisms:
\begin{itemize}
	\item \emph{Morphisms} in $\LLL(F)$ between two objects are given by premorphisms between those two objects modulo homotopy.  
	\item The identity morphism of an object $(x,A)$ in $\LLL(F)$ is given by $[(\id_A, \id_{x}, \id_A)_A]$.
	\item Given two morphisms 
	\begin{equation*}\rule{\leftmargin}{0in} 
		\begin{split}
		[(u_1, f_1, v_1)_{C_1}]:&(x,A) \rightarrow (y,B),\\
		[(u_2,f_2,v_2)_{C_2}]: &(y, B) \rightarrow (z,C),
		\end{split}
	\end{equation*} 
	the composite $[(u_2,f_2,v_2)_{C_2}] \circ [(u_1, f_1, v_1)_{C_1}]$ in $\LLL(F)$ is given by the morphism $[(s_1 \circ u_1, f , s_2 \circ v_2)_D]$ for $s_i: C_i \rightarrow D$ for $i=1,2$ such that $s_1 \circ v_1 = s_2 \circ u_2$ and where the morphism $f: F(s_1 \circ u_1)(x) \rightarrow F(s_2 \circ v_2)(z)$ is defined as the following composite:
	\begin{equation*} \rule{\leftmargin}{0in}
	\begin{aligned}
		&F(s_1 u_1)(x) \cong F(s_1) (F(u_1) (x)) \xrightarrow{F(s_1)(f_1)} F(s_1) (F(v_1)(y)) \cong F(s_1 v_1)(y) = \\
		&F(s_2 u_2)(y) \cong F(s_2) (F(u_2)(y)) \xrightarrow{F(s_2) (f_2)} F(s_2)  (F(v_2) (z)) \cong F(s_2  v_2) (z).
	\end{aligned}
	\end{equation*}
	One can check that this is well-defined. We do not write the details, but essentially, the argument goes as follows. By chosing two different possible representatives of the composition, one can find natural candidates for a homotopy between them by using the fact that $\aaa$ is filtered. In order to check that any of these natural choices is indeed a homotopy, one just needs to use the fact that the isomorphisms $F(g) \circ F(f) \Rightarrow F(g\circ f)$ are natural in both $f$ and $g$ for $f$ and $g$ composable morphisms in $\iii$.  
\end{itemize}

The category $\LLL(F)$ fulfills the universal property (\ref{univprop}) above together with the pseudonatural transformation $\Lambda: F \Rightarrow \LLL(F)$ given by:
\begin{itemize}
	\item for each $A \in \iii$, the functor $\Lambda_A: F(A) \to \LLL(F)$ sends an object $x 
	\in F(A)$ to $(x,A) \in \LLL(F)$ and a morphism $f:x \to y$ in $F(A)$ to $[(1_A,f,1_A)_A]$;
	\item for each $w: A \to B$ in $\iii$, the invertible 2-morphism $\Lambda_w: \Lambda_B \circ F(w) \Rightarrow \Lambda_A$ is given, for each $x \in F(A)$, by the isomorphism
	\begin{equation*}
		[(1_B,1_{F(w)(x)},w)_B]: (F(w)(x),B) \to (x,A)
	\end{equation*}  
	whose inverse is $[(w,1_{F(w)(x)},1_B)_B]$.
\end{itemize}

%Observe that, when $F$ is a strict $2$-functor, this provides a different construction for the usual filtered colimit of categories as in \cite[Example 5.2.2.f]{handbook-categorical-algebra2}.  \doubt{It is not completely clear, as the universal properties are not the same! (look at the strictification)}

We are interested in $k$-linear categories and their $2$-filtered bicolimits. More explicitely, we are interested in 2-filtered bicolimits of pseudofunctors $F: \iii \to \Cat$ that take values in the 2-category $\Cat(k)$ of small $k$-linear categories with $k$-linear functors and $k$-linear natural transformations, that is, in functors $F: \iii \to \Cat$ that factor through the forgetful functor $\Cat(k) \to \Cat$. In general, the 2-filtered bicolimit $\LLL(F)$ under these hypothesis will not necessarily be $k$-linear, as it is also the case for the classical Grothendieck construction (we point the reader to \cite{hochschild-cohomology-presheaves-map-graded-categories} for an account on a linearized Grothendieck construction). However, for nice choices of the indexing $2$-category $\iii$ and the functor $F$, this will hold true. In particular, it is true for our case of interest:
\begin{proposition}\label{klinear}
Let $\iii$ be a filtered category. Take $F: \iii \rightarrow \Cat$ a pseudofunctor that factors through the forgetful functor $\Cat(k) \to \Cat$. Then, $\LLL(F)$ is a $k$-linear category. 
\begin{proof}
	Consider $(x,A), (y,B) \in \LLL(F)$. The class $\LLL(F)((x,A), (y,B))$ has a natural structure of $k$-module induced from the $k$-linear structure of the values of $F$. Indeed, given two morphisms
	\begin{equation*}
			[(u_1, f_1, v_1)_{C_1}], [(u_2,f_2,v_2)_{C_2}]: (x,A) \rightarrow (y,B)
	\end{equation*}
	and an element $\lambda \in k$, we define 
	\begin{equation*}
		[(u_1, f_1, v_1)_{C_1}] + \lambda [(u_2,f_2,v_2)_{C_2}] = [(w_1 \circ u_1, F(w_1)(f_1) + \lambda F(w_2)(f_2), w_2 \circ v_2)_{C}]
	\end{equation*}
	where $w_1: C_1 \rightarrow C$, $w_2: C_2 \rightarrow C$ and $w_1 \circ u_1 = w_2 \circ u_2$ and $w_1 \circ v_1 = w_2 \circ v_2$. Observe that such $w_1, w_2$ exist because $\iii$ is a filtered category. An easy check shows that this is well-defined and does not depend on the choice of $w_1$ and $w_2$. The fact that it provides a $k$-module structure is directly deduced from the fact that, by hypothesis, for each object $C\in \iii$, $F(C)$ is a $k$-linear category and for each morphism $D \rightarrow E$ in $\iii$, the functor $F(D \rightarrow E): F(D) \rightarrow F(E)$ is $k$-linear.
	
	In addition, one can easily show that the composition is $k$-linear. To show this, consider morphisms 
	\begin{equation*}
		\begin{aligned}
			[(u_1, f_1, v_1)_{C_1}], [(u_2,f_2,v_2)_{C_2}]: (x,A) &\rightarrow (y,B),\\
			[(s_1, g_1, t_1)_{D_1}], [(s_2,g_2,t_2)_{D_2}]: (y,B) &\rightarrow (z,C)
		\end{aligned}
	\end{equation*}
	and elements $\lambda,\lambda' \in k$. Without loss of generality we can assume that 
	\begin{equation*}
		\begin{aligned}
			D &= C_1 = C_2 = D_1 = D_2,\\
			u &= u_1 = u_2: A \rightarrow D,\\
			v &= v_1 = v_2 = s_1 = s_2: B \rightarrow D\\
			\text{and}\\
			t&=t_1=t_2: C \rightarrow D.
		\end{aligned}
	\end{equation*}
	We have that 
	\begin{equation*}
		\begin{aligned}
			[(v, g_1, t)_D] \circ & \left( [(u, f_1, v)_D] +\lambda [(u,f_2,v)_D]\right) \\
			&= [(v, g_1, t)_D] \circ [(u, f_1 + \lambda f_2, v)_D]\\
			&= [(u, g_1 \circ (f_1 + \lambda f_2),t)_D]\\
			&= [(u, g_1 \circ f_1 + \lambda g_1 \circ f_2,t)_D]\\
			&= [(v, g_1, t)_D] \circ [(u,f_1,v)_D] + \lambda [(v, g_1, t)_D] \circ [(u,f_2,v)_D].
		\end{aligned}
	\end{equation*} 
	Similarly, one proves that $$\left( [(v, g_1, t)_D] + \lambda' [(v, g_2, t]_D \right) \circ [(u, f_1, v)_D] = [(v, g_1, t)_D].$$
	Hence, the composition is $k$-linear as desired. 
	\end{proof}
\end{proposition}

\begin{remark}
	Observe that one could define a \emph{$k$-linear 2-filtered bicolimit} by replacing in \Cref{defbicolimit} above $\Cat$ by $\Cat(k)$ and the category $\Psnat(F, \ccc)$ by its $k$-linear analogue. Notice then that, given a pseudofunctor $F$ as in \Cref{klinear} above, we have that the ($2$-)filtered bicolimit $\LLL(F)$ coincides with the $k$-linear ($2$-)filtered bicolimit of $F$. In other words, the forgetful functor $U: \Cat(k) \to \Cat$ preserves and reflects filtered bicolimits, as it happens with the forgetful functor $\Ab \to \Set$ \cite[Prop 2.13.5]{handbook-categorical-algebra1}.
\end{remark}

\section{Locally presentable linear categories and Grothendieck categories as filtered bicolimits of small categories}\label{pargrothasbicolim}
In this section we show, based on \S\ref{parfilteredbicolim}, that every locally presentable $k$-linear category can be written as a filtered bicolimit of small categories, and hence in particular this holds for any Grothendieck $k$-linear category. Furthermore, we show that Grothendieck $k$-linear categories can be written as $k$-linear filtered bicolimits of certain filtered categories of $k$-linear sites. 

\begin{theorem}\label{filteredsubcategories}
	Let $\ccc$ be a category which is a union of full small subcategories indexed by a directed poset. Then, $\ccc$ is the filtered bicolimit of that family of subcategories. More precisely,  
	if $\ccc$ is a category such that $\ccc \cong \bigcup_{i \in I} \ccc_i$, where $\ccc_i \subseteq \ccc$ are full small subcategories, $I$ is a directed poset and $\ccc_i \subseteq \ccc_j$ if and only if $i \leq j$, then $\ccc$ is a filtered bicolimit of the family $(\ccc_i)_{i \in I}$.
	\begin{proof}
		Denote by $\iii$ the filtered category given by the directed poset $I$, and denote by $\iota_{i,j}: \ccc_i \hookrightarrow \ccc_j$ the natural embeddings for $i \leq j$.
		We define the $2$-functor 
		\begin{equation}
		F_{\ccc}: \iii \rightarrow \Cat
		\end{equation}
		to be given by
		$F_{\ccc}(i) = \ccc_i$ for every $i \in \iii$ and $F_{\ccc}(i \leq j)= \iota_{i,j}: \ccc_i \hookrightarrow \ccc_{j}$ for every morphism $i \leq j$ in $\iii$.
		
		We build a functor
		$$\phi: \LLL(F_{\ccc}) \to \ccc$$
		defined as follows: 
		\begin{itemize}
			\item $\phi(x, i) = x \in \ccc_i \subseteq \ccc$ for every $(x,i) \in \LLL(F_{\ccc})$;
			\item $\phi\left( [(i\leq k, f: \iota_{i,k}(x) \to \iota_{j,k}(y),j\leq k)] \right)  = f \in \ccc_k(\iota_{i,k}(x), \iota_{j,k}(y)) = \ccc(x, y)$ for every morphism $[(i\leq k, f: \iota_{i,k}(x) \to \iota_{j,k}(y),j\leq k)]:(x,i) \to (y, j)$ in $\LLL(F_{\ccc})$.
		\end{itemize} 
		One can readily check that this is well-defined and it defines a functor. As $\ccc \cong \bigcup_{i \in \iii} \ccc_i$, one trivially has that this functor is essentially surjective, and because $\iota_{i,j}$ are fully faithful for all $i, j \in \iii$, one easily deduces that the functor is fully faithful. 
	\end{proof}
\end{theorem}

\begin{remark}\label{linearbicolimit}
	Assume $\ccc$ is $k$-linear, and hence $\ccc_i$ is $k$-linear for all $i \in \iii$ and so are the fully faithful functors $\iota_{i,j}: \ccc_i \subseteq \ccc_j$ for $i \leq j$ in $\iii$. Then we have that both $\iii$ and $F_{\ccc}$ are in the hypothesis of Proposition \ref{klinear} above, hence $\LLL(F_{\ccc})$ is $k$-linear. Observe that the functor $\phi$ defined in the proof above is as well $k$-linear, and thus $\ccc$ and $\LLL(F_{\ccc})$ are equivalent as $k$-linear categories. 
\end{remark}

\begin{remark}
	\Cref{filteredsubcategories} is just a bicategorical parallel to the easy $1$-categorical easy statement that says that every object in a category that can be written as a union of a filtered system of subobjects is in particular the filtered colimit of these system of subobjects. 
\end{remark}

%As indicated above, we denote by $\regcard$ the category obtained from the totally ordered set of regular cardinals. Recall that given $\ccc$ a locally presentable $k$-linear category and $\alpha \in \regcard$, we denote by $\ccc^{\alpha}$ the full $k$-linear subcategory of $\ccc$ consisting of the $\alpha$-presentable objects. Recall that $\ccc^{\alpha}$ is an $\alpha$-cocomplete category and that if $\alpha \leq \beta$, we have a fully faithful embedding $\ccc^{\alpha} \subseteq \ccc^{\beta}$.  For these and other basic facts concerning locally presentable categories we point the reader to \cite[Ch 1]{locally-presentable-accessible-categories}.

\begin{corollary}\label{thmgrothbicolim}
	Let $\ccc$ be a locally presentable $k$-linear category. Then $\ccc$ is the $k$-linear filtered bicolimit of its family of subcategories of locally presentable objects $(\ccc^{\alpha})_{\alpha}$, where $\alpha$ varies in $\regcard$. 
	\begin{proof}
	As $\ccc \cong \bigcup_{\alpha \in \regcard} \ccc^\alpha$ (see \Cref{remark:raising-cardinality}), the result follows from \Cref{filteredsubcategories} and \Cref{linearbicolimit}.
	\end{proof}
\end{corollary}
\begin{remark}
	Observe that \Cref{thmgrothbicolim} is true for Grothendieck k-linear categories, as
	they are an in particular locally presentable k-linear categories.
\end{remark}

We now introduce another presentation of a Grothendieck $k$-linear category as a filtered bicolimit, where the indexing filtered category will be a certain category of $k$-linear sites. 

Consider the category $\jjj_{\ccc}$ defined as follows:
\begin{itemize}
	\item Objects of $\jjj_{\ccc}$ are given by 
	$$\{u_{\AAA}: (\AAA, \calT_{\AAA}) \to \ccc \,\,|\,\, (\AAA, \calT_{\AAA}) \,\,\text{small }k\text{-linear site}, u_{\AAA} \in \LC\}$$ 
	where $\ccc$ is endowed with the canonical topology (see \Cref{canonicaltop}). For readibility, we will frequently omit the topology from our notation and write $u_{\AAA}: \AAA \to \ccc$.
	\item Morphisms between two objects $u_{\AAA}: \AAA \to \ccc$ and $u_{\BBB}: \BBB \to \ccc$ are given by the $k$-linear functors $f: \AAA \to \BBB$ which belong to $\LC$ and such that $u_{\BBB} \circ f = u_{\AAA}$. We write $f:u_{\AAA} \to u_{\BBB}$.
\end{itemize}
One can readily check this is a well-defined category as a direct consequence of \Cref{propclassLC}.
\begin{proposition}\label{jcfiltered}
	Given a Grothendieck category $\ccc$, the category $\jjj_{\ccc}$ constructed above is filtered. 
	\begin{proof}
	Observe that the category $\jjj_{\ccc}$ is not empty. Given two objects $u_{\AAA}: \AAA \to \ccc$ and $u_{\BBB}: \BBB \to \ccc$, we want to find a third object $u_{\CCC}: \CCC \to \ccc$ and morphisms $f: u_{\AAA} \to u_{\CCC}$ and $g: u_{\BBB} \to u_{\CCC}$. Let $\CCC$ be the full subcategory of $\ccc$ with objects $\{u_{\AAA}(A)\}_{A \in \AAA} \cup \{u_{\BBB}(B)\}_{B \in \BBB}$ and endow it with the topology induced by the canonical topology in $\ccc$. We hence have that the embedding $u_{\CCC}: \CCC \to \ccc$ is an LC morphism. Now consider the corestrictions $u_{\AAA}: \AAA \to \CCC$ and $u_{\BBB}: \BBB \to \CCC$. We trivially have that these define morphisms $u_{\AAA} \to u_{\CCC}$ and $u_{\BBB} \to u_{\CCC}$.
	
	Consider now two morphisms $f,g: u_{\AAA} \to u_{\BBB}$. We want to find an object $u_{\CCC}: \CCC \to \ccc$ and a morphism $h: u_{\BBB} \to u_{\CCC}$ equalizing $f$ and $g$. Take $\CCC$ to be the full subcategory of $\ccc$ with objects $u_{\BBB}(\BBB)$ and endow it with the (restriction of the) canonical topology. Take $u_{\CCC}: \CCC \to \ccc$ to be the embedding (which is LC) and $h: \BBB \to \CCC$ the corestriction of $u_{\BBB}: \BBB \to \ccc$ to $\CCC$. Then, by definition, one has that $h u_{\CCC} = u_{\BBB}$. Furthermore, we have that $h f = h g$ as a direct consequence of the fact that $u_{\BBB} f = u_{\AAA} = u_{\BBB} g$.
	
	We can thus conclude that $\jjj_{\ccc}$ is a filtered category. 
	\end{proof}
\end{proposition}

We now consider the $k$-linear functor $G_{\ccc}: \jjj_{\ccc} \to \Cat(k)$ given by forgetting the ``slice structure'', i.e. defined by sending each object $u_{\AAA}: \AAA \to \ccc$ to the small $k$-linear category $\AAA$ and each morphism $f: u_{\AAA} \to u_{\BBB}$ to itself seen as a $k$-linear morphism $f: \AAA \to \BBB$. 

We proceed to describe $\LLL( G_{\ccc})$ using the construction from \S\ref{parfilteredbicolim} above. Observe that the description in this case will be simplified because $G_{\ccc}$ is a strict functor:
\begin{itemize}
	\item Objects are $\{(x,u_{\AAA}: \AAA \to \ccc) \,\,|\,\, (u_{\AAA}: \AAA \to \ccc) \in \jjj_{\ccc}, x \in \AAA \}$ 
	\item Morphisms $(x,u_{\AAA}: \AAA \to \ccc) \to (y,u_{\BBB}: \BBB \to \ccc)$ are given by homotopy classes of triples $(u, f, v)$ where $u: u_{\AAA} \to u_{\CCC}$, $v: u_{\BBB} \to u_{\CCC}$ are morphisms in $\jjj_{\ccc}$ and $f: u(x) \to v(y)$ is a morphism in $\CCC$. As before, we use the notation $(u, f, v)_{u_{\CCC}}$ to make explicit the codomain of $u$ and $v$. Two morphisms $(u_1, f, v_1)_{u_{\CCC_1}}$ and $(u_2, g, v_2)_{u_{\CCC_2}}$ are homotopic if there exist morphisms $w_1: u_{\CCC_1} \to u_{\CCC}$ and $w_2: u_{\CCC_2} \to u_{\CCC}$ such that $w_1 u_1 = w_2 u_2$, $w_1 v_1 = w_2 v_2$ and $w_1(f)= w_2(g)$. As in \S\ref{parfilteredbicolim}, we denote the homotopy class of $(u_1, f, v_1)_{u_{\CCC_1}}$ by $[(u_1, f, v_1)_{u_{\CCC_1}}]$.
\end{itemize}

\begin{theorem}\label{thmbicolimitsites}
	Given a Grothendieck category $\ccc$, we have that $\ccc$ is the $k$-linear filtered bicolimit of $G_{\ccc}$.
	\begin{proof}
		Observe that $G_{\ccc}$ factors through $\Cat(k)$ and hence, by \Cref{klinear}, we have that $\LLL(G_{\ccc})$ is a $k$-linear category. To conclude, it suffices to construct a $k$-linear equivalence $\psi_{\ccc}: \LLL (G_{\ccc}) \to \ccc$.
		With the notation introduced in \S\ref{parfilteredbicolim} for the objects and morphisms of $\LLL (G_{\ccc})$, we consider the following assignations: 
		\begin{itemize}
			\item To every object $(x, u_{\AAA}: \AAA \to \ccc) \in \jjj_{\ccc}$ we assign the object $$\psi_{\ccc}(x, u_{\AAA}: \AAA \to \ccc) \coloneqq u_{\AAA}(x) \in \ccc.$$
			\item To every morphism $[(u,f,v)_{u_{\CCC}}]:(x,u_{\AAA}:\AAA \to \ccc) \to (y,u_{\BBB}: \BBB \to \ccc)$ we assign the morphism $$\psi_{\ccc}([(u,f,v)_{u_{\CCC}}]) \coloneqq \left( u_{\CCC}(f): u_{\AAA}(x) = u_{\CCC} u(x)  \to u_{\CCC} v(y) = u_{\BBB}(y) \right) .$$
		\end{itemize}
		We have the following:
		\begin{enumerate}
			\item The assignation on morphisms is well-defined. Consider two homotopic morphisms $(u_1, f_1, v_1)_{u_{\CCC_1}}$, $(u_2, f_2, v_2)_{u_{\CCC_2}} :(x,\AAA\to \ccc) \to (y,\BBB \to \ccc)$, this is, there exist morphisms $w_1: u_{\CCC_1} \to u_{\CCC}$ and $w_2: u_{\CCC_2} \to u_{\CCC}$ such that $w_1 u_1 = w_2 u_2$, $w_1 v_1 = w_2 v_2$ and $w_1(f_1) = w_2(f_2)$. We want to show that $u_{\CCC_1}(f_1) = u_{\CCC_2}(f_2)$. Observe that $u_{\CCC_1} (f_1) = u_{\CCC} w_1(f_1) = u_{\CCC} w_2(f_2) = u_{\CCC_2}(f_2)$, as desired.
			\item The assignations define a functor $\psi_{\ccc}: \LLL(G_{\ccc}) \to \ccc$. First observe that it preserves identities. Indeed, the identity morphism $[(\id_{\AAA}, \id_{x}, \id_{\AAA})_{u_{\AAA}}]$ of the object $(x, u_{\AAA}: \AAA \to \ccc)$ gets sent to $u_{\AAA}(\id_{x}) = \id_{u_{\AAA}(x)}$. 
			We now check that it preserves compositions. Consider two composable morphisms $[ (u_1,f_1,v_1)_{u_{\CCC_1}} ]: (\AAA,u_{\AAA}) \to (\BBB,u_{\BBB})$ and $[(u_2,f_2,v_2)_{u_{\CCC_2}}]: (\BBB,u_{\BBB}) \to (\CCC,u_{\CCC})$ in $\LLL(G_{\ccc})$. Because $\jjj_{\ccc}$ is filtered, we can assume that $u_{\DDD} = u_{\CCC_1} = u_{\CCC_2}$, i.e. that the morphisms $u_1, v_1, u_2, v_2$ in $\jjj_{\ccc}$ have the same codomain $u_{\DDD}: \DDD \to \ccc$, and that $v_1 = u_2$.
			Their composite in $\LLL(G_{\ccc})$ is given by $[(u_1,f_2 \circ f_1,v_2)_{u_{\DDD}}]$ and gets sent to $u_{\DDD} (f_2 f_1)$. 
			On the other hand, we have that $[ (u_1,f_1,v_1)_{u_{\DDD}} ]$ and $[(u_2,f_2,v_2)_{u_{\DDD}}]$ get sent to $u_{\DDD}(f_1)$ and $u_{\DDD}(f_2)$ respectively, whose composite is $u_{\DDD}(f_2) u_{\DDD}(f_1)$. As $u_{\DDD}:\DDD \to \ccc$ is a functor, we have that $u_{\DDD}(f_2 \circ f_1) = u_{\DDD}(f_2) u_{\DDD}(f_1)$ as desired.
			\item The functor $\psi_{\ccc}$ is $k$-linear. Let $[ (u_1,f_1,v_1)_{u_{\CCC_1}} ], [(u_2,f_2,v_2)_{u_{\CCC_2}}]: (\AAA,u_{\AAA}) \to (\BBB,u_{\BBB})$ be two morphisms in $\LLL(G_{\ccc})$. Because $\jjj_{\ccc}$ is filtered, we may assume that $u_{\CCC} = u_{\CCC_1} = u_{\CCC_2}$ and that $u = u_1 = u_2$, $v = v_1 = v_2$. Then, given $\lambda \in k$, we have that 
			\begin{equation*}\rule{\leftmargin}{0in}
				\psi_{\ccc}([ (u,f_1,v)_{u_{\DDD}} ] + \lambda [(u,f_2,v)_{u_{\DDD}}]) = \psi_{\ccc}([u, f_1 + \lambda f_2,v]_{u_{\DDD}}) = u_{\DDD}(f_1 + \lambda f_2).
			\end{equation*}
			On the other hand, we have that 
			\begin{equation*}\rule{\leftmargin}{0in}
				\psi_{\ccc}([ (u,f_1,v)_{u_{\DDD}} ]) + \lambda \psi_{\ccc}([(u,f_2,v)_{u_{\DDD}}]) = u_{\DDD}(f_1) + \lambda u_{\DDD}(f_2).
			\end{equation*}
			As $u_{\DDD}$ is a $k$-linear functor, we have that $u_{\DDD}(f_1 + \lambda f_2) = u_{\DDD}(f_1) + \lambda u_{\DDD}(f_2)$ as desired.
			\item The functor $\psi_{\ccc}$ is essentially surjective. Let $y$ be an object in $\ccc$. Consider the small full subcategory $\AAA$ of $\ccc$ with objects $\{y\} \cup \{g\}_{g \in G}$, where $G$ is a small set of generators of $\ccc$. We endow $\AAA$ with the topology induced by the canonical topology in $\ccc$. Then the embedding $\iota: \AAA \to \ccc$ is trivially an LC morphism and hence we have that $\psi_{\ccc}(y,\iota:\AAA \to \ccc) = \iota(y) = y$, as desired.
			\item The functor $\psi_{\ccc}$ is faithful. Let $[(u_1,f_1,v_1)_{u_{
			\CCC_1}}]$, $[(u_2,f_2,v_2)_{u_{\CCC_2}}]: (x,u_{\AAA}) \to (y,u_{\BBB})$ be two morphisms, such that $u_{\CCC_1}(f_1) = u_{\CCC_2}(f_2)$. Consider the full subcategory $\DDD$ of $\ccc$ with objects $\{u_{\CCC_1}(\CCC_1)\} \cup \{u_{\CCC_2}(\CCC_2)\}$ endowed with the topology induced by the canonical topology in $\ccc$ and the associated embedding $\iota: \DDD \to \ccc$. Then we have that $u_{\CCC_1}$ and $u_{\CCC_2}$ factor through $\iota$:
			\begin{equation*}\rule{\leftmargin}{0in}
				\begin{tikzcd}
				\CCC_1 \arrow[rd, "\tilde{u}_{\CCC_1}"] \arrow[rrd, "u_{\CCC_1}", bend left] &  &  \\
				& \DDD \arrow[r, "\iota"] & \ccc. \\
				\CCC_2 \arrow[ru, "\tilde{u}_{\CCC_2}"'] \arrow[rru, "u_{\CCC_2}"', bend right] &  & 
				\end{tikzcd}
			\end{equation*}
			Observe that, because $u_{\CCC_1}(f_1) = u_{\CCC_2}(f_2)$ and $\iota$ is an embedding, we have that $\tilde{u}_{\CCC_1} u_1 = \tilde{u}_{\CCC_2}  u_2$, $\tilde{u}_{\CCC_1}v_1 = \tilde{u}_{\CCC_2}v_2$ and $\tilde{u}_{\CCC_1}(f_1) = \tilde{u}_{\CCC_2}(f_2)$. This implies that $(u_1,f_1,v_1)_{u_{\CCC_1}}$ and $(u_2,f_2,v_2)_{u_{\CCC_2}}$ are homotopic, and hence $[(u_1,f_1,v_1)_{u_{\CCC_1}}] = [(u_2,f_2,v_2)_{u_{\CCC_2}}]$.  
			\item The functor $\psi_{\ccc}$ is full. Consider $(x, u_{\AAA}: \AAA \to \ccc), (y, u_{\BBB}: \BBB \to \ccc) \in \jjj_{\ccc}$ and a morphism $f: u_{\AAA}(x) \to u_{\BBB}(y)$ in $\ccc$. As previously, consider the full subcategory $\CCC$ of $\ccc$ with objects $\{u_{\AAA}(\AAA)\} \cup \{u_{\BBB}(\BBB)\}$ endowed with the topology induced by the canonical topology in $\ccc$ and the embedding $\iota: \CCC \to \ccc$, which is an LC morphism. Then, we have that $u_{\AAA}$ and $u_{\BBB}$ factor through $\iota$ as above:
			\begin{equation*}\rule{\leftmargin}{0in}
			\begin{tikzcd}
			\AAA \arrow[rd, "\tilde{u}_{\AAA}"] \arrow[rrd, "u_{\AAA}", bend left] &  &  \\
			& \CCC \arrow[r, "\iota"] & \ccc. \\
			\BBB \arrow[ru, "\tilde{u}_{\BBB}"'] \arrow[rru, "u_{\BBB}"', bend right] &  & 
			\end{tikzcd}
			\end{equation*}
			We can hence consider $\tilde{f}: \tilde{u}_{\AAA}(x) \to \tilde{u}_{\BBB}(y)$ the image of $f$ via the isomorphism $\ccc(\iota \tilde{u}_{\AAA} (x) , \iota \tilde{u}_{\BBB}(y)) \cong \CCC(\tilde{u}_{\AAA}(x), \tilde{u}_{\BBB}(y))$ induced by $\iota$. Take the morphism $[(\tilde{u}_{\AAA}, \tilde{f}, \tilde{u}_{\BBB})_{\iota = u_{\CCC}}]: (x, u_{\AAA}:\AAA \to \ccc) \to (y,u_{\BBB}:\BBB \to \ccc )$. By construction we have that $\psi_{\ccc}$ sends this morphism to $\iota (\tilde{f}) = f$, proving fulness. 
		\end{enumerate}
	We hence have proven that $\psi_{\ccc}:\LLL(G_{\ccc}) \to \ccc$ is a $k$-linear equivalence of categories as desired.
	\end{proof}
\end{theorem}

\begin{remark}
Observe that, in order to recover any locally presentable category $\ccc$ as a filtered bicolimit using the construction from \Cref{thmgrothbicolim}, we can always use the same filtered category, namely the category $\regcard$ associated to the total ordered class of small regular cardinals. Notice that this is not the case for this last presentation of Grothendieck categories provided by \Cref{thmbicolimitsites}, as the filtered category $\jjj_{\ccc}$ depends on the Grothendieck category $\ccc$ we want to recover.
\end{remark}

\section{The tensor product of locally presentable linear categories as a filtered bicolimit}\label{tensorproduct} 
In this section we analyse the tensor product of locally presentable $k$-linear categories in terms of the realization of locally presentable categories as filtered bicolimits provided in \S\ref{pargrothasbicolim}. In particular, this applies to the tensor product of Grothendieck $k$-linear categories.
 
First observe that, as a direct consequence of \Cref{filteredsubcategories}, we have that $\aaa \boxtimes_{\mathsf{LP}} \bbb$ is the filtered bicolimit of the family $((\aaa \boxtimes_{\mathsf{LP}} \bbb)^{\alpha})_{\alpha \in \regcard}$ and that from $\alpha \geq \kappa$ we have that $(\aaa \boxtimes_{\mathsf{LP}} \bbb)^{\alpha} \cong \aaa^{\alpha} \otimes_{\alpha} \bbb^{\alpha}$. However, it is not directly obvious whether the fully faithful functors $\aaa^{\alpha} \otimes_{\alpha} \bbb^{\alpha} \hookrightarrow \aaa^{\beta} \otimes_{\beta} \bbb^{\beta}$ with $\beta \geq \alpha \geq \kappa$ provided by the inclusion $(\aaa \boxtimes_{\mathsf{LP}} \bbb)^{\alpha} \subseteq (\aaa \boxtimes_{\mathsf{LP}} \bbb)^{\beta}$ coincide with the canonical functors $f_{\alpha,\beta}$ from \eqref{transitionfunctors}, and hence whether the filtered bicolimit is compatible with the $\alpha$-cocomplete tensor products for $\alpha$ varying in $\regcard$. This is precisely the content of the following theorem.

\begin{theorem}\label{remfullyfaithful}
	Let $\aaa$, $\bbb$ be two locally presentable $k$-linear categories and choose the smallest regular cardinal $\kappa$ such that both $\aaa$ and $\bbb$ are locally $\kappa$-presentable.
	Then, for all $\alpha, \beta \in \regcard$ such that $\beta \geq \alpha \geq \kappa$ the canonical functor $$f_{\alpha, \beta}: \aaa^{\alpha} \otimes_\alpha \bbb^{\alpha} \rightarrow \aaa^{\beta} \otimes_\beta \bbb^{\beta}$$ defined in \eqref{transitionfunctors} is fully faithful. In particular, the functor $f_{\alpha,\beta}$ coincides, up to the equivalences $(\aaa \boxtimes_{\LP} \bbb)^{\alpha} \cong \aaa^{\alpha} \otimes_{\alpha} \bbb^{\alpha}$ and $(\aaa \boxtimes_{\LP} \bbb)^{\beta} \cong \aaa^{\beta} \otimes_{\beta} \bbb^{\beta}$, with the canonical inclusion $\iota_{\alpha,\beta}^{\aaa \boxtimes_{\LP} \bbb}:(\aaa \boxtimes_{\LP} \bbb)^{\alpha} \subseteq (\aaa \boxtimes_{\LP} \bbb)^{\beta}$.
	
	\begin{proof}
		Let $\iota \coloneqq \iota_{\alpha,\beta}^{\aaa} \otimes \iota_{\alpha,\beta}^{\bbb}$ as in \Cref{proposition:2variablesdiagramlocallypresentable}. First observe that in the diagram
		\begin{equation*}
			\begin{tikzcd}
				\aaa^\alpha \otimes \bbb^\alpha \arrow[dd, "{R_{\alpha,\alpha} Y_{\alpha,\alpha}}" description] \arrow[rrrr, "\iota" description] \arrow[rd, "{u_{\aaa^\alpha,\bbb^\alpha}}" description] &&&& \aaa^\beta \otimes \bbb^\beta \arrow[dd,"{R_{\beta,\beta} Y_{\beta,\beta}}" description] \arrow[ld, "{u_{\aaa^\beta,\bbb^\beta}}" description] \\
				& \aaa^\alpha \otimes_\alpha \bbb^\alpha \arrow[rr, "{f_{\alpha,\beta}}" description] \arrow[ld, "k_\alpha" description, hook] && \aaa^\beta \otimes_\beta \bbb^\beta \arrow[rd, "k_\beta" description, hook] &\\
				{\Lex_{\alpha,\alpha}(\aaa^\alpha,\bbb^\alpha)} \arrow[rrrr, "\iota^s" description] &&&& {\Lex_{\beta,\beta}(\aaa^\beta,\bbb^\beta)}            
			\end{tikzcd}
		\end{equation*}
		the outer square commutes by \Cref{proposition:2variablesdiagramlocallypresentable}, the left and the right triangles commute by \Cref{not:alpha-alpha-vs-alpha} and the upper square commutes by \eqref{eq:raising-cardinality-Kelly}. Therefore, we have that
		\begin{equation*}
			k_\beta f_{\alpha,\beta} u_{\aaa^\alpha,\bbb^\alpha} = \iota^s k_\alpha u_{\aaa^\alpha,\bbb^\alpha}.
		\end{equation*}
		As $k_\beta, f_{\alpha,\beta}, k_\alpha, \iota^s$ all preserve $\alpha$-small colimits and every element in $\aaa^\alpha \otimes_{\alpha} \bbb^\alpha$ is an $\alpha$-small colimit of elements in the image of $u_{\aaa^\alpha,\bbb^\alpha}$ by \Cref{prop:explicit-construction-Kelytp}, we conclude that the lower square also commutes, namely that
		\begin{equation*}
			k_\beta f_{\alpha,\beta} = \iota^s k_\alpha.
		\end{equation*}
		Thus, because $ k_\alpha, k_\beta$ and $\iota^s$ are fully faithful, we can conclude that $f_{\alpha,\beta}$ is also fully faithful as desired. 
		
		We know we have an equivalence $E:\aaa \boxtimes_{\LP} \bbb \to \Lex_{\alpha,\alpha}(\aaa^\alpha,\bbb^\alpha)$, and define $F: \aaa \boxtimes_{\Groth} \bbb \to \Lex_{\beta,\beta}(\aaa^\beta,\bbb^\beta) $ to be the equivalence given by $F \coloneqq \iota^s E$. Then, we have a commutative diagram
		\begin{equation*}
			\begin{tikzcd}[column sep=small]
				& \aaa^\alpha \otimes_\alpha \bbb^\alpha \arrow[rr, "{f_{\alpha,\beta}}" description, hook] \arrow[dd, "k_\alpha" description, hook] &{}                                                                                                                                                              & \aaa^\beta \otimes_\beta \bbb^\beta \arrow[dd, "k_\beta" description, hook] &                                                                                                                                                                                       \\
				&                                         &                                                                                                                                                                 &                                                                             &                                                                                                                                                                                       \\
				& {\Lex_{\alpha,\alpha}(\aaa^\alpha,\bbb^\alpha)} \arrow[rr, "\iota^s" description]                                                
				& {} \arrow[uu, phantom, "\circlearrowleft" description]                                                                                                                   & {\Lex_{\beta,\beta}(\aaa^\beta,\bbb^\beta)}                                 &                                                                                                                                                                                       \\
				(\aaa \boxtimes_{\LP} \bbb)^\alpha \arrow[ruuu, "E'" description, bend left] \arrow[rr, "\iota^{\aaa \boxtimes_{\LP} \bbb}_\alpha" description, hook] \arrow[rrrr, "{\iota^{\aaa \boxtimes_{\LP} \bbb}_{\alpha,\beta}}" description, bend right, hook] \arrow[ru, phantom, "\circlearrowleft" description] &                                                                                                                              & \aaa \boxtimes_{\LP} \bbb \arrow[lu, "E" description] \arrow[ru, "F" description] \arrow[u, phantom, "\circlearrowleft" description] \arrow[dd, phantom, "\circlearrowleft" description] &                                                                             
				& (\aaa \boxtimes_{\LP} \bbb)^\beta, \arrow[ll, "\iota^{\aaa \boxtimes_{\LP} \bbb}_{\beta}" description, hook] \arrow[luuu, "F'" description, bend right] \arrow[lu, phantom, "\circlearrowleft" description] 
				\\
				&&&&
				\\
				&&{}&&                                                                    
			\end{tikzcd}
		\end{equation*}
		where $E'$ (resp. $F'$) is the equivalence obtained by restricting $E$ (resp. $F$) to $(\aaa \boxtimes_{\Groth} \bbb)^\alpha$ (resp. $(\aaa \boxtimes_{\Groth} \bbb)^\beta$). Consequently, we have that
		\begin{equation*}
			k_\beta f_{\alpha,\beta} E' = k_\beta F' \iota_{\alpha,\beta}^{\aaa \boxtimes_{\Groth} \bbb},
		\end{equation*}
		and as $k_\beta$ is injective on objects and fully faithful, we can conclude that
		\begin{equation*}
			f_{\alpha,\beta} E' = F' \iota_{\alpha,\beta}^{\aaa \boxtimes_{\Groth} \bbb},
		\end{equation*}
		as desired.
	\end{proof}
\end{theorem}

We can now define 
\begin{equation}\label{pseudofunctorcd}
F_{\aaa,\bbb}: \regcard \rightarrow \Cat 
\end{equation}
the pseudofunctor given by 
\begin{itemize}
	\item $F_{\aaa,\bbb}(\alpha)= F_{\aaa}(\alpha) \otimes_{\alpha} F_{\bbb}(\alpha) = \aaa^{\alpha} \otimes_{\alpha} \bbb^{\alpha}$ for every $\alpha \in \regcard$;
	\item $F_{\aaa,\bbb}(\alpha \leq \beta)= f_{\alpha, \beta}$ for every morphism $\alpha \leq \beta$ in $\regcard$;
\end{itemize}
with the notation from \S\ref{pargrothasbicolim} above.
\begin{remark}
	Observe that $F_{\aaa,\bbb}$ is a pseudofunctor and not a strict $2$-functor.
\end{remark}
%\begin{remark}\label{remfullyfaithful}
%	Observe that, as the canonical functor $\aaa^{\alpha} \otimes \bbb^{\alpha} \rightarrow \aaa^{\alpha} \otimes_{\alpha} \bbb^{\alpha}$ is fully faithful and $i^{\aaa}_{\alpha, \beta}, i^{\bbb}_{\alpha, \beta}$ are fully faithful as well (and hence so is $i^{\aaa}_{\alpha, \beta} \otimes i^{\bbb}_{\alpha, \beta}$), we have that $i^{\aaa}_{\alpha, \beta} \otimes_{\alpha} i^{\bbb}_{\alpha, \beta}$ is also fully faithful.
%\end{remark}
\begin{theorem}\label{productasbicolim}
	Given $\aaa$ and $\bbb$ two locally presentable $k$-linear categories, one has that $\aaa \boxtimes_{\LP} \bbb$ is the $k$-linear filtered bicolimit of $F_{\aaa, \bbb}$.
	\begin{proof}
		By Proposition \ref{klinear}, $\LLL (F_{\aaa,\bbb})$ is a $k$-linear category. In addition, we know that there exists a smallest $\kappa \in \regcard$ such that $\aaa \boxtimes_{\LP} \bbb$ is locally $\alpha$-presentable and there is an equivalence $(\aaa \boxtimes_{\LP} \bbb)^{\alpha} \cong \aaa^{\alpha} \otimes_{\alpha} \bbb^{\alpha}$ for every $\alpha \geq \kappa$, which from this point on, and for simplicity, we will treat as an identity. Consider $\kappa$ the smallest regular cardinal with such property. We build a $k$-linear functor
		$$\phi: \LLL(F_{\aaa, \bbb}) \to \aaa \boxtimes_{\LP} \bbb$$
	as follows.
	
	For an object $(x,\alpha) \in \LLL(F_{\aaa,\bbb})$, we put 
	$$\phi(x, \alpha) = \begin{cases} 
	x \in \aaa^{\alpha} \otimes_{\alpha} \bbb^{\alpha} \subseteq \aaa \boxtimes_{\LP} \bbb &\text{if } \alpha \geq \kappa \\
	f_{\alpha,\kappa} (x) \in \aaa_{\kappa} \otimes_{\kappa} \bbb_{\kappa} \subseteq \aaa \boxtimes_{\LP} \bbb &\text{if } \alpha < \kappa
	\end{cases}$$
		
	For a morphism $[(\alpha \leq \gamma,g,\beta \leq \gamma)]: (x,\alpha) \to (y, \beta)$, making use of \Cref{remfullyfaithful}, we define $\phi\left( [(\alpha \leq \gamma,g,\beta \leq \gamma)]\right) $ as follows:	   
	\begin{itemize}
		\item If $\alpha, \beta \geq \kappa$, we put
		\begin{equation*}
		\begin{aligned}
		\phi\left( [(\alpha \leq \gamma,g,\beta \leq \gamma)]\right)  = g &\in \aaa_{\gamma} \otimes_{\gamma} \bbb_{\gamma}(f_{\alpha,\gamma}(x),f_{\beta,\gamma}(y) ) \\
		&\cong \aaa \boxtimes_{\LP} \bbb(x,y).
		\end{aligned}
		\end{equation*}
	 	Observe that in this case $\gamma \geq \kappa$ holds.
		\item If $\alpha < \kappa$ and $\beta \geq \kappa$, we put
		\begin{equation*}
		\begin{aligned}
		\phi\left( [(\alpha \leq \gamma,g,\beta \leq \gamma)]\right)  = g &\in \aaa_{\gamma} \otimes_{\gamma} \bbb_{\gamma}(f_{\alpha,\gamma}(x),f_{\beta,\gamma}(y) )\\ 
		&\cong \aaa \boxtimes_{\LP} \bbb(f_{\alpha,\kappa}(x),y).
		\end{aligned}
		\end{equation*} 
		Observe that in this case $\gamma \geq \kappa$ holds.
		\item If $\alpha \geq \kappa$ and $\beta < \kappa$, we put 
		\begin{equation*}
		\begin{aligned}
		\phi\left( [(\alpha \leq \gamma,g,\beta \leq \gamma)]\right) = g &\in \aaa_{\gamma} \otimes_{\gamma} \bbb_{\gamma}(f_{\alpha,\gamma}(x),f_{\beta,\gamma}(y) )\\
		&\cong \aaa \boxtimes_{\LP} \bbb(x,f_{\beta,\kappa}(y)).
		\end{aligned}
		\end{equation*}
		Observe that in this case $\gamma \geq \kappa$ holds.
		\item If $\alpha, \beta < \kappa$ and $\gamma \geq \kappa$, we put
		\begin{equation*}
			\begin{aligned}
				\phi\left( [(\alpha \leq \gamma,g,\beta \leq \gamma)]\right)  = g &\in \aaa_{\gamma} \otimes_{\gamma} \bbb_{\gamma}(f_{\alpha,\gamma}(x),f_{\beta,\gamma}(y) )\\
				&\cong \aaa \boxtimes_{\LP} \bbb(f_{\alpha,\kappa}(x),f_{\beta,\kappa}(y)).
				\end{aligned}
			\end{equation*}
		\item If $\alpha, \beta, \gamma < \kappa$, we put:
		 \begin{equation*}
		 	\begin{aligned}
		 		\phi\left( [(\alpha \leq \gamma,g,\beta \leq \gamma)]\right)  = f_{\gamma,\kappa}(g) &\in \aaa_{\kappa} \otimes_{\kappa} \bbb_{\kappa}(f_{\alpha,\kappa}(x),f_{\beta,\kappa}(y) )\\ 
		 		&\cong \aaa \boxtimes_{\LP} \bbb(f_{\alpha,\kappa}(x),f_{\beta,\kappa}(y)).
		 		\end{aligned}
	 		\end{equation*}
	\end{itemize}
Observe this functor is well-defined and $k$-linear. In addition, we have that 
$$\aaa \boxtimes_{\LP} \bbb \cong \bigcup_{\alpha \in \regcard} (\aaa \boxtimes_{\LP} \bbb)^{\alpha} \cong \bigcup_{\alpha \geq \kappa} (\aaa \boxtimes_{\LP} \bbb)^{\alpha} \cong \bigcup_{\alpha \geq \kappa} \aaa^{\alpha} \otimes_{\alpha} \bbb^{\alpha}.$$ 
Consequently, the functor is essentially surjective. We also have that all the transition functors $f_{\alpha,\beta}$ are fully-faithful for $\beta \geq \alpha \geq \kappa$ by \Cref{remfullyfaithful} above, hence one can conclude that the functor is fully-faithful as desired. 
	\end{proof}
\end{theorem}

\begin{remark}\label{remark:tplinearsites-vs-tpGrothendieck}
	One may wonder if an analogous approach would allow to obtain a realization of the tensor product of Grothendieck categories as a filtered bicolimit by using, instead of the $\alpha$-Kelly tensor product, the tensor product $\boxtimes_{\sites}$ of linear sites and LC morphisms from \S\ref{linearsites} and, instead of the realization of Grothendieck categories as filtered bicolimits of $\alpha$-presentable objects from \Cref{thmgrothbicolim}, the realization of Grothendieck categories as filtered bicolimits of linear sites from \Cref{thmbicolimitsites}. We will explain why this is not the case. Roughly, the argument goes as follows: 
	
	Let $\aaa,\bbb$ be Grothendieck $k$-linear categories. We use the notation introduced in \S\ref{pargrothasbicolim} for the rest of the remark. Consider the filtered categories $\jjj_{\aaa}$ (resp. $\jjj_{\bbb}$) with objects the LC morphisms $u: (\CCC, \calT_{\CCC}) \to (\aaa,\calT_{\aaa,\mathrm{can}})$, (resp. the LC morphisms $v:(\DDD, \calT_{\DDD}) \to (\bbb,\calT_{\bbb,\mathrm{can}})$). Denote by $\jjj_{\aaa,\bbb}$ the category with objects given by tensor products $u \otimes v: (\CCC \otimes \DDD, \calT_{\CCC} \boxtimes_{\sites} \calT_{\DDD}) \to (\aaa \otimes \bbb, \calT_{\aaa,\mathrm{can}} \boxtimes_{\sites} \calT_{\bbb,\mathrm{can}})$ of objects in $\jjj_{\aaa}$ with objects of $\jjj_{\bbb}$, and with morphisms given by the tensor product of LC morhisms. In particular, we have that $\jjj_{\aaa,\bbb}$ is filtered. Moreover, as $\LC$ is closed under composition (see \Cref{propclassLC}) and tensor products (see \Cref{tpLC}), composition with the LC morphism 
	$$s_{\aaa,\bbb}:(\aaa \otimes \bbb, \calT_{\aaa,\mathrm{can}} \boxtimes_{\sites} \calT_{\bbb,\mathrm{can}}) \to (\aaa \boxtimes_{\Groth} \bbb, \calT_{\aaa \boxtimes_{\Groth} \bbb,\mathrm{can}})$$ 
	defines a faithful functor $T: \jjj_{\aaa,\bbb} \to \jjj_{\aaa \boxtimes_{\Groth} \bbb}$, where $\jjj_{\aaa \boxtimes_{\Groth} \bbb}$ is the category with objects the LC morphisms $(\EEE,\calT_{\EEE}) \to (\aaa \boxtimes_{\Groth} \bbb, \calT_{\aaa \boxtimes_{\Groth} \bbb,\mathrm{can}})$.  Consider now the functor $\bar{G}_{\aaa,\bbb}$ given by the composite
	\begin{equation*}
		\begin{tikzcd}
			{\jjj_{\aaa,\bbb}} \arrow[r, "T"] & \jjj_{\aaa \boxtimes_{\Groth} \bbb} \arrow[r, "G_{\aaa \boxtimes_{\Groth} \bbb}"] & \Cat.
		\end{tikzcd}
	\end{equation*}
	One can then construct a $k$-linear functor 
	$$\bar{\psi}_{\aaa,\bbb}: \LLL (\bar{G}_{\aaa,\bbb}) \to \aaa \boxtimes_{\Groth} \bbb,$$  
	by restricting the functor $\psi_{\aaa \boxtimes_{\Groth} \bbb}: \LLL (G_{\aaa \boxtimes_{\Groth} \bbb}) \to \aaa \boxtimes_{\Groth} \bbb$ constructed in the proof of \cref{thmbicolimitsites} above.
	However, $\bar{\psi}_{\aaa,\bbb}$ is not an equivalence. This follows from the fact that the natural functor $\aaa \otimes \bbb \to \aaa \boxtimes_{\Groth} \bbb$ is not essentially surjective (consider, for example, the natural functor $\Mod(\AAA) \otimes \Mod(\BBB) \to \Mod(\AAA \otimes \BBB)$, which is not essentially surjective in general\footnote{A counterexample of the essential surjectivity of this functor can be found, for example, in the answer provided by Pierre-Yves Galliard in the following entry of Mathematics StackExchange: \url{https://math.stackexchange.com/questions/104222/modules-over-a-tensor-product}}), which implies that $\bar{\psi}_{\aaa,\bbb}$ is also not essentially surjective. 
\end{remark}

\section{The functoriality, associativity and symmetry of the tensor product of locally presentable linear categories in terms of filtered bicolimits}\label{parfunctassocsym}
In this section we describe the functoriality, associativity and symmetry of the tensor product of locally presentable linear categories (and thus in particular that of Grothendieck categories) with respect to cocontinuous functors in terms of the functoriality, associativity and symmetry of the $\alpha$-Kelly tensor products, making use of \Cref{productasbicolim} above. In particular, this provides an advantage when computing the tensor product of cocontinuous morphisms when we have control of the subcategories of presentable objects, as illustrated by \Cref{example:tp-qcoh} below. 

\begin{definition}\label{def:tp-cocontinuous-functors}
	Consider locally presentable $k$-linear categories $\aaa, \bbb, \ccc$ and $\ddd$ and cocontinuous functors $F: \aaa \to \ccc$ and $F': \bbb \to \ddd$. We define 
	\begin{equation*}
		F \boxtimes_{\LP} F' : \aaa \boxtimes_{\LP} \bbb \to \ccc \boxtimes_\LP \ddd
	\end{equation*}
	to be the cocontinuous functor obtained from applying the universal property of $\boxtimes_\LC$ from \eqref{eq:universal-property-tplp} to the (cococontinuous in each variable) functor
	\begin{equation*}
		\begin{tikzcd}
			\aaa \otimes \bbb \arrow[r,"F \otimes F'"] &\ccc \otimes \ddd \arrow[r] &\ccc \boxtimes_\LP \ddd.
		\end{tikzcd}
	\end{equation*}
\end{definition}

In what follows, we are going to provide a description of the tensor product $\boxtimes_{\LP}$ of cocontinuous linear functors making use of the filtered bicolimit presentation of $\boxtimes_{\LP}$ from \Cref{productasbicolim}. In particular, we will use the functoriality of the $\alpha$-Kelly tensor product:
\begin{definition}\label{def:tp-alpha-cocontinuous-functors}
	Consider $\alpha$-cocomplete $k$-linear categories $\AAA, \BBB, \CCC$ and $\DDD$ and $\alpha$-cocontinuous functors $f: \AAA \to \CCC$ and $f': \BBB \to \DDD$. We define 
	\begin{equation*}
		f \otimes_{\alpha} f' : \AAA \otimes_{\alpha} \BBB \to \CCC \otimes_\alpha \DDD
	\end{equation*}
	to be the $\alpha$-cocontinuous functor obtained from applying the universal property of $\otimes_\alpha$ from \eqref{alphaprod} to the ($\alpha$-cococontinuous in each variable) functor
	\begin{equation*}
		\begin{tikzcd}
			\AAA \otimes \BBB \arrow[r,"f \otimes f'"] &\CCC \otimes \DDD \arrow[r] &\CCC \otimes_\alpha \DDD.
		\end{tikzcd}
	\end{equation*}
\end{definition}

Consider $\aaa$, $\bbb$ two locally presentable categories and a regular cardinal $\alpha$. Recall that a linear functor $F: \aaa \to \bbb$ is said to \emph{have rank $\alpha$} if it preserves $\alpha$-filtered colimits \cite[\S 5.5]{handbook-categorical-algebra2}. It is trivial to see that if a functor $F$ has rank $\alpha$, then it has rank $\beta$ for every $\beta \geq \alpha$. We say a functor \emph{has rank} if there exists a regular cardinal $\alpha$ such that it has rank $\alpha$. 
We have the following useful proposition.
\begin{proposition}[{\cite[Prop. 5.5.6]{handbook-categorical-algebra2}}]\label{leftadjointrank}
	Let $G: \aaa \to \bbb$ a functor between locally presentable $k$-linear categories. If $G$ has a left adjoint, then $G$ has a rank.
\end{proposition}
The following is easy to show, but we provide a proof for the convenience of the reader.
\begin{proposition}
	Let $F: \aaa \to \bbb$ a cocontinuous $k$-linear functor between locally presentable $k$-linear categories. Then, there exists a regular cardinal $\alpha$ such that $F(\aaa^{\beta}) \subseteq \bbb^{\beta}$ for every $\beta \geq \alpha$.
	\begin{proof}
		By the dual of the Special Adjoint Functor Theorem \cite[Thm. 3.3.4]{handbook-categorical-algebra1}, we have that $F$ has a right adjoint $G$. In particular, by \Cref{leftadjointrank}, $G$ has rank. Fix the smallest $\alpha$ such that $G$ has rank $\alpha$.  Then, given an element $C \in \aaa^{\alpha}$, we have that 
		\begin{equation*}
			\begin{aligned}
				\bbb(F(C), \colim_i D_i)  &= \aaa(C, G(\colim_i D_i))\\
				&= \aaa(C, \colim_i G(D_i) )\\
				&= \colim_i \aaa (C, G(D_i))\\
				&= \colim_i (F(C), D_i)
			\end{aligned}
	\end{equation*}
		where $\colim_i D_i$ is any $\alpha$-filtered colimit in $\bbb$. Hence $F(C) \in \bbb^{\alpha}$ as desired.
	\end{proof}
\end{proposition}
 
\begin{remark}
	Given $F: \aaa \to \bbb$ as in the proposition, note that the restriction-corestriction $F_{\beta}: \aaa^{\beta} \rightarrow \bbb^{\beta}$ of $F$ is $\beta$-cocontinuous for all $\beta \geq \alpha$.
\end{remark} 

We first show that the functoriality of the tensor product $\boxtimes_{\LP}$ and that of the tensor products $\otimes_{\alpha}$ are nicely compatible in the following sense:
\begin{proposition}\label{prop:compatibility-functoriality}
	Consider locally presentable $k$-linear categories $\aaa, \bbb, \ccc$ and $\ddd$ and cocontinuous $k$-linear functors $F: \aaa \to \ccc$ and $F': \bbb \to \ddd$. Let $\kappa$ be the smallest regular cardinal for which both $F$ and $F'$ preserve $\alpha$-presentable objects for every $\alpha \geq \kappa$.
	For every $\gamma \geq \beta \geq \kappa$, the diagrams
	\begin{equation}\label{eq:functoriality-compatibility-cardinality}
		\begin{tikzcd}
			\aaa^\beta \otimes_\beta \bbb^\beta \arrow[rr, "F_\beta \otimes_\beta F'_\beta" description] \arrow[d, "{f_{\beta,\gamma}}"', hook] &  & \ccc^\beta \otimes_\beta \bbb^\beta \arrow[d, "{f_{\beta,\gamma}}", hook] \\
			\aaa^\gamma \otimes_\gamma \bbb^\gamma \arrow[rr, "F_\gamma \otimes_\gamma F'_\gamma" description]                                                     &  & \ccc^\gamma \otimes_\gamma \ddd^\gamma                                                   
		\end{tikzcd}
	\end{equation}
	and
	\begin{equation}\label{eq:functoriality-compatibility-twotp}
		\begin{tikzcd}
			\aaa^\beta \otimes_\beta \bbb^\beta \arrow[rrr, "F_{\beta} \otimes_{\beta} F'_{\beta}" description] \arrow[d, hook] &  &  & \ccc^\beta \otimes_{\beta} \ddd^\beta \arrow[d, hook] \\
			\aaa\boxtimes_\LP \bbb \arrow[rrr, "F \boxtimes_\LP F'" description]                                                                  &  &  & \ccc \boxtimes_\LP \ddd                                 
		\end{tikzcd}
	\end{equation}
	are commutative.
	\begin{proof}
		Consider the diagram
		\begin{equation*}
			\begin{tikzcd}
				\aaa^\beta \otimes \bbb^\beta \arrow[rd, "{u_{\aaa^\beta,\bbb^\beta}}" description] \arrow[rrrrr, "F_{\beta} \otimes F_{\beta}" description] \arrow[ddd, "{\iota^\aaa_{\beta,\gamma} \otimes \iota^\bbb_{\beta,\gamma}}" description] &                                                                                                                                                              &  &  &                                                                                             & \ccc^\beta \otimes \ddd^\beta \arrow[ld, "{u_{\ccc^\beta,\ddd^\beta}}" description] \arrow[ddd, "{\iota^\ccc_{\beta,\gamma} \otimes \iota^\ddd_{\beta,\gamma}}" description] \\
				& \aaa^\beta \otimes_\beta \bbb^\beta \arrow[rrr, "F_{\beta} \otimes_{\beta} F'_{\beta}" description] \arrow[d, "{f_{\beta,\gamma}}"', hook] &  &  & \ccc^\beta \otimes_{\beta} \ddd^\beta \arrow[d, "{f_{\beta,\gamma}}", hook] &                                                                                                                                                                                    \\
				& \aaa^\gamma \otimes_\gamma \bbb^\gamma \arrow[rrr, "F_\gamma \otimes_\gamma F'_\gamma" description]                                                          &  &  & \ccc^\gamma \otimes_\gamma \ddd^\gamma                                                      &                                                                                                                                                                                    \\
				\aaa^\gamma \otimes \bbb^\gamma \arrow[ru, "{u_{\aaa^\gamma,\bbb^\gamma}}" description] \arrow[rrrrr, "F_\gamma \otimes F'_\gamma" description]                                                                                      &                                                                                                                                                              &  &  &                                                                                             & \ccc^\gamma \otimes \ddd^\gamma \arrow[lu, "{u_{\ccc^\gamma,\ddd^\gamma}}" description]                                                                                           
			\end{tikzcd}
		\end{equation*}
		We want to show that the centre square, which is precisely \eqref{eq:functoriality-compatibility-cardinality}, is commutative. Observe that the upper and lower square are commutative by definition, as so is the outer square. Moreover, the left and right squares are commutative by virtue of \eqref{eq:raising-cardinality-Kelly}. Therefore, as all the functors in the centre square preserve $\beta$-small colimits and $\aaa^\beta \otimes_\beta \bbb^\beta$ is generated under $\beta$-small colimits by elements in the image of $u_{\aaa^\beta,\bbb^\beta}$, we can conclude that the centre diagram is also commutative, as desired. 
		
		Similarly, in order to prove \eqref{eq:functoriality-compatibility-twotp}, one can consider the diagram 
		\begin{equation*}
			\begin{tikzcd}
				\aaa^\beta \otimes \bbb^\beta \arrow[rd, "{u_{\aaa^\beta,\bbb^\beta}}" description] \arrow[rrrrr, "F_{\beta} \otimes F_{\beta}" description] \arrow[ddd, "\iota^\aaa_\beta \otimes \iota^\bbb_\beta" description] &                                                                                                                                       &  &  &                                                          & \ccc^\beta \otimes \ddd^\beta \arrow[ld, "{u_{\ccc^\beta,\ddd^\beta}}" description] \arrow[ddd, "\iota^\ccc_\beta \otimes \iota^\ddd_\beta" description] \\
				& \aaa^\beta \otimes_\beta \bbb^\beta \arrow[rrr, "F_{\beta} \otimes_{\beta} F'_{\beta}" description] \arrow[d, hook] &  &  & \ccc^\beta \otimes_{\beta} \ddd^\beta \arrow[d, hook] &                                                                                                                                                                \\
				& \aaa\boxtimes_\LP \bbb \arrow[rrr, "F \boxtimes_\LP F'" description]                                                                  &  &  & \ccc \boxtimes_\LP \ddd                                  &                                                                                                                                                                \\
				\aaa \otimes \bbb \arrow[ru, "{s_{\aaa,\bbb}}" description] \arrow[rrrrr, "F \otimes F'" description]                                                                                                                                 &                                                                                                                                       &  &  &                                                          & \ccc \otimes \ddd \arrow[lu, "{s_{\ccc,\ddd}}" description]                                                                                                   
			\end{tikzcd}
		\end{equation*}
	The upper, lower and outer square are commutative by definition. In addition, the left and right squares are commutative by virtue of \Cref{rem:compatibility-universal-properties}. Therefore, as all the functors in the centre square preserve $\beta$-small colimits and $\aaa^\beta \otimes_\beta \bbb^\beta$ is generated under $\beta$-small colimits by elements in the image of $u_{\aaa^\beta,\bbb^\beta}$, we can conclude that the centre diagram, which is precisely \eqref{eq:functoriality-compatibility-twotp}, is also commutative, as we wanted to show. 
	\end{proof}
\end{proposition}

Let $F: \aaa \to \ccc$ and $F': \bbb \to \ddd$ be as in \Cref{prop:compatibility-functoriality}. We define a $k$-linear pseudonatural transformation 
$$\Phi:F_{\aaa, \bbb} \Rightarrow \ccc \boxtimes_{\LP} \ddd,$$
where $F_{\aaa,\bbb}$ is defined as in \eqref{pseudofunctorcd}, as follows. For each $\alpha \in \regcard$, we put
\begin{itemize}
	\item If $\alpha < \kappa$, we define $\Phi_{\alpha}$ as the natural composite $$\aaa^{\alpha} \otimes_{\alpha} \bbb^{\alpha} \xrightarrow{f_{\alpha,\kappa}} \aaa^{\kappa} \otimes_{\kappa} \bbb^{\kappa} \xrightarrow{F_{\kappa} \otimes_{\kappa} F'_{\kappa}} \ccc^{\kappa} \otimes_{\kappa} \ddd^{\kappa} \hookrightarrow \ccc \boxtimes_{\LP} \ddd;$$
	\item If $\alpha \geq \kappa$, we define $\Phi_{\alpha}$ as the natural composite $$\aaa^{\alpha} \otimes_{\alpha} \bbb^{\alpha} \xrightarrow{F_{\alpha} \otimes_{\alpha} F'_{\alpha}} \ccc^{\alpha} \otimes_{\alpha} \ddd^{\alpha} \hookrightarrow \ccc \boxtimes_{\LP} \ddd;$$
\end{itemize}
where, for any regular cardinal $\gamma$, $F_{\gamma} \otimes_{\gamma} F'_{\gamma}: \aaa^{\gamma} \otimes_{\gamma} \bbb^{\gamma} \to  \ccc^{\gamma} \otimes_{\gamma} \ddd^{\gamma}$ denotes the natural functor obtained from $F_{\gamma} \otimes F'_{\gamma}: \aaa^{\gamma} \otimes \bbb^{\gamma} \to  \ccc^{\gamma} \otimes \ddd^{\gamma}$ via the universal property of $\otimes_{\gamma}$.
For each morphism $\alpha \leq \beta$ in $\regcard$, we set the invertible natural transformations 
$$\Phi_{\alpha \leq \beta}: \Phi_{\beta} \circ F_{\aaa, \bbb}(\alpha \leq \beta) \Rightarrow \Phi_{\alpha}$$ 
to be the corresponding identities (observe that both sides of the natural transformation coincide as a direct consequence of \Cref{prop:compatibility-functoriality}). 

This construction provides a way to construct the tensor product $F \boxtimes_\LP F'$. Indeed, we have the following: 
\begin{proposition}\label{functoriality}
	Let $F: \aaa \to \ccc$ and $F': \bbb \to \ddd$ be cocontinuous $k$-linear functors between locally presentable $k$-linear categories. The functor associated to the pseudonatural transformation $\Phi: F_{\aaa, \bbb} \Rightarrow \ccc \boxtimes_{\LP} \ddd$ above via the universal property of $\LLL(F_{\aaa,\bbb})$ \eqref{univprop} coincides with $F \boxtimes_\LP F'$.
	\begin{proof}
		By virtue of \eqref{univprop}, in order to conclude it is enough to show that the pseudonatural transformation $\Phi: F \Rightarrow \ccc \boxtimes_\LP \ddd$ is isomorphic to the pseudonatural transformation 
		\begin{equation}\label{eq:pseudonatural-transformation-functoriality}
			\begin{tikzcd}
				F_{\aaa,\bbb} \arrow[r, "\Lambda", Rightarrow] &   \LLL(F_{\aaa,\bbb}) \cong \aaa \boxtimes_\LP \bbb \arrow[r, "F\boxtimes_{\LP}F'", Rightarrow] & \ccc \boxtimes_{\LP} \ddd,
			\end{tikzcd}
		\end{equation} 
		where $\Lambda$ is the universal pseudonatural transformation associated to the identity of $\LLL(F_{\aaa,\bbb})$ (see \S\ref{parfilteredbicolim}) and we are considering the categories $\aaa \boxtimes_{\LP} \bbb$ and $\ccc \boxtimes_\LP \ddd$ as constant $2$-functors from $\regcard$ to $\Cat(k)$ and $F \boxtimes_\LP F'$ as the corresponding constant pseudonatural transformation. 
		
		Let's denote the pseudonatural transformation \eqref{eq:pseudonatural-transformation-functoriality} by $\Psi$. We are going to show that $\Psi$ and $\Phi$ are not only isomorphic, but they are actually the same pseudonatural transformation. Let $\kappa_1$ be the smallest regular cardinal such that $\aaa$ and $\bbb$ are locally $\kappa_1$-presentable and $\kappa_2 \geq \kappa_1$ the smallest regular cardinal such that $\aaa,\bbb,\ccc$ and $\ddd$ are locally $\kappa_2$ presentable and $F,F'$ have rank $\kappa_2$. 
		
		A straightforward computation shows that, given $\alpha \in \regcard$, the $k$-linear functor $\Psi_\alpha: \aaa^\alpha \otimes_{\alpha} \bbb^\alpha \to \ccc \boxtimes_{\LP} \ddd$ is given by:
		\begin{itemize}
			\item if $\alpha \geq \kappa_1$, $\Psi_{\alpha} = (F\boxtimes_{\LP}F')|_{\aaa^{\alpha} \otimes_{\alpha} \bbb^{\alpha}}$;
			\item if $\alpha < \kappa_1$, $\Psi_{\alpha} = (F\boxtimes_{\LP}F')|_{\aaa^{\kappa_1} \otimes_{\kappa_1} \bbb^{\kappa_1}} \circ f_{\alpha,\kappa_1}$.   
		\end{itemize} 
	
		Moreover, given $\alpha \leq \beta$, one can easily check that the invertible $2$-morphism $\Psi_{\alpha \leq \beta}: \Psi_\beta \circ F_{\aaa,\bbb}(\alpha \leq \beta) \Rightarrow \Psi_\alpha$ is given by:
		\begin{itemize}
			\item if $\beta \geq \alpha \geq \kappa_1$, $\Psi_{\alpha \leq \beta}: (F\boxtimes_{\LP}F')|_{\aaa^{\beta} \otimes_{\beta} \bbb^\beta} \circ f_{\alpha,\beta} \Rightarrow (F\boxtimes_{\LP}F')|_{\aaa^{\alpha} \otimes_{\alpha} \bbb^\alpha}$ is the identity;
			\item if $\beta \geq \kappa_1 > \alpha$, $\Psi_{\alpha \leq \beta} : (F\boxtimes_{\LP}F')|_{\aaa^{\beta} \otimes_{\beta} \bbb^\beta} \circ f_{\alpha,\beta} \Rightarrow (F\boxtimes_{\LP}F')|_{\aaa^{\kappa_1} \otimes_{\kappa_1} \bbb^{\kappa_1}} \circ f_{\alpha,\kappa_1}$ is the identity (observe that $f_{\alpha,\beta} = f_{\kappa_1,\beta} \circ f_{\alpha,\beta}$);
			\item if $\kappa_1 > \beta \geq \alpha$: $\Psi_{\alpha \leq \beta} : (F\boxtimes_{\LP}F')|_{\aaa^{\kappa_1} \otimes_{\kappa_1} \bbb^{\kappa_1}} \circ f_{\beta,\kappa_1} \circ f_{\alpha,\beta} \Rightarrow (F\boxtimes_{\LP}F')|_{\aaa^{\kappa_1} \otimes_{\kappa_1} \bbb^{\kappa_1}} \circ f_{\alpha,\kappa_1}$ is the identity as well.   
		\end{itemize} 
		Then, using \Cref{prop:compatibility-functoriality} one can readily check that $\Psi$ coincides with $\Phi$, as desired.
	\end{proof}
\end{proposition}

%\begin{remark}
%	Note that $F \boxtimes_{\Groth} F'$ is also cocontinuous. The filtered nature of the bicolimit plays an important role in the proof. Roughly it can be shown as follows. Consider $\colim_i X_i$ the colimit of a small family of objects in $\aaa \boxtimes_{\Groth} \bbb$. Then, we can choose a regular cardinal $\alpha$ such that $\colim_i X_i$ is an $\alpha$-small colimit, all the $X_i$ are $\alpha$-presentable and $F$ and $F'$ preserve $\alpha$-presentable objects. Then we can see $\colim_i X_i$ as an element in $\aaa^{\alpha} \otimes_{\alpha} \bbb^{\alpha}$ and we have that $$F\boxtimes_{\Groth} F' (\colim_i X_i) = F_{\alpha} \otimes_{\alpha} F_{\alpha}' (\colim_i X_i) = \colim_i (F_{\alpha} \otimes_{\alpha} F_{\alpha}')(X_i) = \colim_i (F\boxtimes_{\Groth} F') (X_i),$$
%	where we have used that $F_{\alpha} \otimes_{\alpha} F'_{\alpha}: \aaa^{\alpha} \otimes_{\alpha} \bbb^{\alpha} \to \ccc^{\alpha} \otimes_{\alpha} \ddd^{\alpha}$ preserves $\alpha$-small colimits by the universal property of $\otimes_{\alpha}$.
%\end{remark}

The description of the tensor product of cocontinuous functors in terms of the Kelly tensor product of presentable objects provides a useful tool in order to perform computations when one has control of the presentable objects. Let us illustrate this with the following geometrical example.
\begin{example}\label{example:tp-qcoh}
Let $f_1: X_1 \to Y_1$ be a morphism of $k$-schemes where $X_1$ and $Y_1$ are quasi-compact and quasi-separated. In particular, we know that:
\begin{itemize}
	\item $\Qch(X_1)$ and $\Qch(Y_1)$ are Grothendieck $k$-linear categories \cite[Tag 077P]{stacks-project};
	\item $\Qch(X_1)$ and $\Qch(Y_1)$ are locally finitely presentable $k$-linear categories \cite[Chap I, Cor 6.9.12]{EGA1};
	\item the pullback functor $f_1^*: \Qch(Y_1) \to \Qch(X_1)$, which is cocontinuous, preserves finitely presentable objects \cite[Chap 0, \S5.2.5]{EGA1}.
\end{itemize}
Let $f_2:X_2 \to Y_2$ be another morphism between quasi-compact quasi-separated $k$-schemes and consider the \emph{external tensor product} $\boxtimes$ of quasi-coherent sheaves \cite[Chap I, Def 3.3.1]{EGA1}, which is given by 
\begin{equation*}
	\begin{aligned}
		\boxtimes: \Qch(X_1) \otimes \Qch(X_2) &\to \Qch(X_1 \times_k X_2)\\
		(F,G) &\mapsto F \boxtimes G \coloneqq p_{X_1}^* (F) \otimes_{\ooo_{X_1 \times_k X_2}} p^*_{X_2}(G), 
	\end{aligned}	
\end{equation*}
where $p_{X_i}: X_1 \times_k X_2 \to X_i$ denote the corresponding projections for $i =1,2$.
From \cite[Thm 6.1]{localizations-tensor-categorires-fiber-products-schemes} we have that the external tensor product $\boxtimes$ exhibits the Grothendieck category $\Qch(X_1 \times_k X_2)$ as the tensor product of Grothendieck $k$-linear categories $\Qch(X_1) \boxtimes_{\Groth} \Qch(X_2)$. In particular, by virtue of \Cref{rem:compatibility-universal-properties}, we have that the restriction-corestriction of the exterior tensor product 
$$\boxtimes: \Qch(X_1)^{\aleph_0} \otimes \Qch(X_2)^{\aleph_0} \to \Qch(X_1 \times_k X_2)^{\aleph_0}$$ 
also exhibits the finitely cocomplete linear category $\Qch(X_1 \times_k X_2)^{\aleph_0}$ as the $\aleph_0$-Kelly tensor product $\Qch(X_1)^{\aleph_0} \otimes_{\aleph_0} \Qch(X_2)^{\aleph_0}$. 

As a consequence of \Cref{prop:compatibility-functoriality}, we have the following commutative diagram:
	\begin{equation*}
		\begin{tikzcd}
			\Qch(Y_1)^{\aleph_0} \otimes \Qch(Y_2)^{\aleph_0} \arrow[d, "\boxtimes" description] \arrow[rrr, "(f^*_1)_{\aleph_0} \otimes (f^*_2)_{\aleph_0}" description] \arrow[dd, "\boxtimes"', bend right=74] &  &  & \Qch(X_1)^{\aleph_0} \otimes \Qch(X_2)^{\aleph_0} \arrow[d, "\boxtimes" description] \arrow[dd, "\boxtimes", bend left=74] \\
			\Qch(Y_1 \times_k Y_2)^{\aleph_0} \arrow[rrr, "(f_1^*)_{\aleph_0} \otimes_{\aleph_0} (f_2^*)_{\aleph_0}" description] \arrow[d, hook]                                       &  &  & \Qch(X_1 \times_k X_2)^{\aleph_0} \arrow[d, hook]                                                                          \\
			\Qch(Y_1 \times_k Y_2) \arrow[rrr, "f_1^* \boxtimes_{\Groth} f_2^*" description]                                                                                                     &  &  & \Qch(X_1 \times_k X_2)                                                                                                    
		\end{tikzcd}
	\end{equation*} 
and thus, for any $F \in \Qch(Y_1)^{\aleph_0}$ and $G \in \Qch(Y_2)^{\aleph_0}$, we have that $$(f_1^* \boxtimes_{\Groth} f_2^*)(F \boxtimes G) = f_1^*(F) \boxtimes f_2^*(G).$$ 
Moreover, notice that this fully determines the cocontinuous functor $f_1^* \boxtimes_{\Groth} f_2^*$. Indeed, every sheaf in $\Qch(Y_1 \times_k Y_2)$ is written as a filtered colimit of sheaves in $\Qch(Y_1 \times_k Y_2)^{\aleph_0}$ and furthermore, every sheaf in $\Qch(Y_1 \times_k Y_2)^{\aleph_0}$ is a finite colimit of sheaves of the form $F \boxtimes G$ where $F \in \Qch(Y_1)^{\aleph_0}$ and $G \in \Qch(Y_2)^{\aleph_0}$. 

This example is particularly interesting in the case in which $Y_1$ and $Y_2$ are locally noetherian, as in that case, for $i =1,2$ we have that $\Qch(Y_i)^{\aleph_0}$ is precisely $\Coh(Y_i)$, the full subcategory of coherent sheaves \cite[Tag 01XZ]{stacks-project}. 
\end{example}

Now, we show how one can recover the associativity and the symmetry of the tensor product of locally presentable linear categories from the same properties of the $\alpha$-Kelly tensor products.

Consider locally presentable $k$-linear categories $\aaa, \bbb$ and $\ccc$.
Define $$F_{(\aaa, \bbb), \ccc}: \regcard \to \Cat$$ 
by $F_{(\aaa, \bbb), \ccc}(\alpha) = (\aaa^{\alpha} \otimes_{\alpha} \bbb^{\alpha}) \otimes_{\alpha} \ccc^{\alpha}$
with the natural transition functors 
$$F_{(\aaa, \bbb), \ccc}(\alpha \leq \beta) :(\aaa^{\alpha} \otimes_{\alpha} \bbb^{\alpha}) \otimes_{\alpha} \ccc^{\alpha} \to (\aaa^{\beta} \otimes_{\beta} \bbb^{\beta}) \otimes_{\beta} \ccc^{\beta}$$
induced by the universal property of $\otimes_{\alpha}$. Analogously, we put 
$$F_{\aaa, (\bbb, \ccc)}: \regcard \to \Cat$$
with $F_{\aaa, (\bbb, \ccc)}(\alpha) = \aaa^{\alpha} \otimes_{\alpha} (\bbb^{\alpha} \otimes_{\alpha} \ccc^{\alpha})$ with the natural transition functors $$F_{\aaa, (\bbb, \ccc)}(\alpha \leq \beta) :\aaa^{\alpha} \otimes_{\alpha} (\bbb^{\alpha} \otimes_{\alpha} \ccc^{\alpha}) \to \aaa^{\beta} \otimes_{\beta} (\bbb^{\beta} \otimes_{\beta} \ccc^{\beta}).$$
In a similar fashion to \Cref{productasbicolim}, one can show that $$\LLL (F_{(\aaa, \bbb), \ccc}) \cong (\aaa \boxtimes_{\LP} \bbb) \boxtimes_{\LP} \ccc,$$ and analogously $$\LLL (F_{\aaa, (\bbb, \ccc)}) \cong \aaa \boxtimes_{\LP} (\bbb \boxtimes_{\LP} \ccc).$$ We know that, for any regular cardinal $\alpha$, the category $\Cat_{\alpha}(k)$ of $\alpha$-cocomplete small categories endowed with $\otimes_{\alpha}$ is a closed monoidal symmetric bicategory. In particular, we have that
$$(\AAA \otimes_{\alpha} \BBB) \otimes_{\alpha} \CCC \cong \AAA \otimes_{\alpha}(\BBB \otimes_{\alpha} \CCC)$$
for all $\alpha$-cocomplete categories $\AAA, \BBB, \CCC$. 

Consequently, there is a canonical isomorphism $F_{(\aaa, \bbb), \ccc}(\alpha) \cong F_{\aaa, (\bbb, \ccc)}(\alpha)$ for each $\alpha$, and it behaves functorially. We thus have
\begin{equation}\label{associativity}
	F_{(\aaa, \bbb), \ccc} \cong F_{\aaa, (\bbb, \ccc)}. 
\end{equation}

\begin{proposition}\label{propassociativity}
	Let $\aaa, \bbb$ and $\ccc$ be locally presentable $k$-linear categories, then, there exists an equivalence
	 $$(\aaa \boxtimes_{\LP} \bbb) \boxtimes_{\LP} \ccc \cong \aaa \boxtimes_{\LP} (\bbb \boxtimes_{\LP} \ccc).$$
	 \begin{proof}
	 	It follows from applying filtered bicolimits to \eqref{associativity}.
	 \end{proof}
\end{proposition}

The argument to deduce the symmetry of the tensor product of locally presentable linear categories from the Kelly tensor products is analogous. 
Consider two locally presentable $k$-linear categories $\aaa$ and $\bbb$. As the monoidal bicategory $(\Cat_{\alpha}(k),\otimes_{\alpha})$ is symmetric, we have
$$\AAA \otimes_{\alpha} \BBB \cong \BBB \otimes_{\alpha}\AAA$$
for all $\alpha$-cocomplete categories $\AAA, \BBB$.

Thus, reasoning as above, we have a canonical isomorphism 
\begin{equation}\label{symmetry}
	F_{\aaa, \bbb} \cong F_{\bbb,\aaa},
\end{equation} 
where $F_{\aaa,\bbb}$ and $F_{\bbb,\aaa}$ are defined as in \eqref{pseudofunctorcd}.
 
\begin{proposition}\label{propsymmetry}
	Let $\aaa, \bbb$ be locally presentable $k$-linear categories. Then, there exists an equivalence
	$$\aaa \boxtimes_\LP \bbb \cong \bbb \boxtimes_\LP \aaa.$$
	 \begin{proof}
		It follows from applying filtered bicolimits to \eqref{symmetry}.
	\end{proof}
\end{proposition}

\printbibliography

\end{document}